\documentclass[11pt]{amsart} 

\usepackage{amssymb,amsmath,amscd,amsfonts,graphicx,latexsym,epsfig}

\usepackage[all]{xy}
\newtheorem{teo}{Theorem}[section]
\newtheorem{lem}[teo]{Lemma}

\newtheorem{cor}[teo]{Corollary}
\newtheorem{exa}[teo]{Example}
\newtheorem{prop}[teo]{Proposition}
\newtheorem{defi}[teo]{Definition}
\newtheorem{ques}[teo]{Question}

\newtheorem{remark}[teo]{Remark}

\newcommand{\Ff}{{\mathcal F}}

\newcommand{\Ll}{{\mathcal L}}
\newcommand{\Mm}{{\mathcal M}}

\newcommand{\Ss}{{\mathcal S}}

\newcommand{\HH}{{\bf H}}

\newcommand{\Z}{{\mathbb Z}}
\newcommand{\R}{{\mathbb R}}

\newcommand{\N}{{\mathbb N}}

\def\Dim{\noindent\emph{Proof. }}
\def\cvd{\hfill$\Box$}
\newcommand{\lra}{\longrightarrow}%{\to}
\newcommand{\B}{\mathsf{B}}
\newcommand{\compl}{\mathsf{C}}
\newcommand{\fr}{\Z^{\ast r}}

\title[]{The topology of Helmholtz domains }

\author[]{R. Benedetti$^1$, R. Frigerio$^1$, R. Ghiloni$^2$}

\address{$^1$ Dipartimento di Matematica \\
Universit\`a di Pisa \\
Largo B.~Pontecorvo 5 \\
56127 Pisa, Italy}

\address{$^2$ Dipartimento di Matematica\\ Universit\`a di Trento\\
Via Sommarive 14\\ 38123 Povo, Italy} 

\email{benedett@dm.unipi.it, frigerio@dm.unipi.it, ghiloni@science.unitn.it}



\begin{document}

\maketitle

\begin{abstract} The goal of this paper is to describe and
clarify as much as possible the 3--dimensional topolo\-gy underlying 
the Helmholtz cuts method, which occurs in a wide theoretic and applied literature about {\it Electromagnetism}, {\it Fluid dynamics} and {\it Elasticity} on domains of the ordinary space $\R^3$. We consider two classes of
bounded domains that satisfy mild boundary conditions and that become
``simple'' after a finite number of disjoint cuts along properly
embedded surfaces. For the first class ({\it Helmholtz}), ``simple''
means that every curl--free smooth vector field admits a potential. For the
second ({\it weakly--Helmholtz}), we only require that a potential
exists for the restriction of every curl--free smooth vector field defined on
the whole initial domain. By means of classical and rather elementary
facts of 3--dimensional geometric and algebraic topology, we give an
exhaustive description of Helmholtz domains, realizing that their
topology is forced to be quite elementary (in particular, Helmholtz
domains with connected boundary are just possibly knotted
handlebodies, and the complement of any non--trivial link is not~Helmholtz).
The discussion about weakly--Helmholtz domains is a bit
more advanced, and their classification appears to be a quite
difficult issue. Nevertheless, we provide several interesting
characterizations of them and, in particular, we point out that the 
class of links with weakly--Helmholtz complements eventually coincides with the
class of the
so--called {\it homology boundary links}, that have been widely studied
in Knot Theory. 
\end{abstract}

\section{Introduction}\label{intro}
{\it Hodge decomposition} is an important analytic structure occurring in a
wide theoretic and applied literature on {\it Electromagnetism}, {\it
Fluid dynamics} and {\it Elasticity} on domains of the ordinary space
$\R^3$ (see a selection of titles in ``Section A'' of our
References). In \cite{CDTG}, one can find a friendly introduction to
this topic. Helmholtz's ``cuts method'' arised in this framework, as
far as we understand, in order to obtain a more effective description
of the Hodge decomposition of the space of $L^2$--vector fields on a
given domain, which could also allow explicit numerical
processings. These ideas can be incorporated in the notion of
so--called {\it Helmholtz domain}. Roughly speaking, a Helmholtz domain
is a bounded
domain that becomes ``simple'' after a finite number of cuts along
disjoint surfaces.  It turns out that there is a bit of indeterminacy
in the literature about the right meaning of ``simple''. Requiring
the domain to be {\it simply connected} certainly suffices. However,
the (possibly weaker) condition consisting in the {\it existence of
potentials for curl--free smooth vector fields} sounds more pertinent
to the actual setting.  Apparently, the relationship between such a
priori different notions is not widely well established. In Section 16
of \cite{CDTG}, one can find a historical account about the way
embryonic forms of homotopy and homology groups of spatial domains had
been introduced by
Helmholtz, Thomson and reconsidered by Maxwell
in the study of Electro and Fluid dynamics. 
Quoting from page 439:

\smallskip

{\it ``Thomson introduced an embryonic version of the one--dimensional
homology $H_1(\Omega)$ in which one countes the number of
``irreconcilable" closed paths inside the domain $\Omega$. This was
subject to the standard confusion of the time between homology and
homotopy of paths: homology was the appropriate notion in this
setting, but the definitions were those of homotopy".}

\smallskip

One could say that such a confusion of the early times somehow
propagated by internal paths till the present days (including true
misunderstandings, see the discussion of Example~\ref{Vou-Bin} below).

\smallskip

On the other hand, spatial domains (whose study includes, for example,
{\it Knot Theory}) re\-present a non--trivial specialization of
3--dimensional manifolds and, since Poincar\'e's {\it Analysis Situs}
(1895) (\cite{Gordon} provides an useful historical account), an
important range of applications of the ideas and techniques of
(3--dimensional) {\it Geometric and Algebraic Topology} developed time by time. 

\smallskip

The first aim of the present largely expository paper is 
to completely clarify the topology of Helmholtz domains, just by applying 
a few classical results or rather elementary facts of 3--dimensional topology.
\smallskip

The first results we recognize (see Theorem \ref{simple-top},
Corollary \ref{cor-simple}) show that, under mild assumptions on the
boundary (e.g. when the boundary is locally Lipschitz,
condition which is usually taken for granted in the literature on
Helmholtz domains), the notions of ``simplicity'' mentioned above are
indeed equivalent to each other. Moreover, it turns out that
simple domains admit a clear and easy description: they are just the
complement of a finite number of disjoint balls in a larger ball. In
the case of polyhedral boundaries, this is due to Borsuk \cite{Borsuk}
(1934). The validity for more general ({\it locally flat}) topological
boundaries depends on later deep results that we will recall in
Theorem \ref{smoothing}. The proof we will provide is based on
elementary properties of the {\it Euler--Poincar\'e characteristic} of
compact surfaces and 3--manifolds and (like in~\cite{Borsuk}) eventually
reduces to the celebrated Alexander Theorem \cite{ALEX} (1924)
asserting that every polyhedral (locally flat indeed) 2--sphere in
$\R^3$ bounds a 3--ball. In \cite{Fox} (1948), Fox obtained Borsuk's
Theorem as a corollary of his {\it reimbedding theorem} (see
Section \ref{reimbedding} below). However, Fox's arguments are
admittedly inspired by 
Alexander's results and techniques.

\smallskip

Once simple domains have been completely described, it is rather
easy to give an exhaustive characterization of general Helmholtz
domains (see Theorem \ref{helm-char}).  In a sense, this is a
disappointing result, as it shows that the topology of Helmholtz
domains is forced to be quite elementary. For
example, Helmholtz domains with connected boundary are just
(possibly knotted) handlebodies, and the complement of any
non--trivial link is not Helmholtz.
\smallskip

In Section~\ref{weakly-Helmholtz}, we introduce and discuss the
strictly larger class of so--called {\it weakly--Helmholtz} domains.
Roughly speaking, such a domain can be cut along a finite number of
disjoint surfaces into subdomains on which curl--free smooth vector fields,
that are defined on the {\it whole} original domain, admit
potentials. We believe that this requirement naturally weakens the Helmholtz
condition, thus allowing to apply the method of cuts to topologically
richer classes of domains. Unlike in the case of 
Helmholtz domains, we are not 
able to give an exhaustive classification of weakly--Helmholtz ones. However, we
will provide several interesting characterizations of weakly--Helmholtz
domains. In particular and remarkably, we realize that the class of
links with weakly--Helmholtz complements eventually coincides with the
class of so--called {\it homology boundary links}. In particular,
every knot and every classical boundary link has weakly--Helmholtz
complement. Homology boundary links are very widely studied in Knot
Theory, and it is a nice occurrence that the Helmholtz cut method
leads to such a distinguished class of links. 

Paper \cite{GHI} is a sort of complement to the present one. It deals
with an effective description of the Hodge decomposition of the space
of $L^2$--vector fields on any bounded domain of $\R^3$ with
sufficiently regular boundary, without making use of any cuts--type
method.

\smallskip

We stress that, from the strict 3--dimensional topology viewpoint, the
results of this paper are largely applications of classical and
well--known facts of Differential/Algebraic/Geo\-metric Topology, that
are usually covered by basic courses on these subjects. This reflects
upon ``Section B'' of our References, that contains well established
books on these subjects, that are exhaustive for our needs. In order
to make the exposition simpler for a reader not too familiar with such
topics, instead of recalling these facts in one comprehensive section,
we have preferred to do it time by time. As already said, the
discussion about Helmholtz domains only needs simple facts about the
Euler--Poincar\'e characteristic (see Section \ref{algtop1:subsec}),
together with Alexander's Theorem.  Very clear and accessible proofs
of this last result are available (e.g. in \cite{HATCH2}).  The discussion
about weakly--Helmholtz domains is a bit more advanced. More
information on the algebraic topology of spatial domains is developed
in Section \ref{algtop2:subsec}, and we will make intensive use of
duality.

On the other hand, we hope that this paper could be of some utility to
people interested in research areas mentioned at the beginning of this
introduction. The r\^ole of
the (algebraic) topology
of domains 
had already been stressed in \cite{CDTG} and \cite{GrKo} (for example in order to
justify the dimension of the Hodge decomposition summands). Hopefully,
the present work should integrate the papers just mentioned, by unfolding the
3--dimensional topology underlying the Helmholtz cuts method.

\smallskip

{\bf Aknowledgements.} The first two authors like to thank R. Ghiloni
for having ``discovered'' the Helmholtz domains literature, and involved
them in the task of clarifying their topology. The third author thanks
especially Alberto Valli, for having introduced him to these themes, and
convinced him of the utility of such a task. He also thanks Ana {Alonso
Rodr{\'\i}guez}, Annaliese Defranceschi, Domenico Luminati, and Valter Moretti 
for helpful conversations.

\section{Domains}\label{domains}
In what follows, {\it smooth} maps (whence, in particular,
diffeomophisms) or manifolds will always assumed to be of class
$C^\infty$.

First a few terminology. The terms ``disk'' and ``ball'' are often used
indifferently, by specifying time by time if they are open or
closed. We prefer here to profit of both terms by stipulating that a
disk is closed and a ball is the open interior of a disk. More
precisely, let $(x_1,x_2,x_3)$ be the usual coordinates of $\R^3$ and
let $D^3$ be the standard $3$--disk $\{(x_1,x_2,x_3) \in \R^3 \, | \,
x_1^2+x_2^2+x_3^2 \leq 1\}$ of~$\R^3$.  Identify $\R^2$ with the plane
$x_3=0$ of $\R^3$ and denote by $D^2$ the standard $2$--disk defined
by $D^2:=D^3 \cap \R^2$.

\begin{defi}
{\rm A subset $X$ of a manifold $M$ homeomorphic to $\R^3$ is a {\it
$($topological$)$ $3$--disk} if, up to homeomorphism, the pair $(M,X)$ is
equivalent to $(\R^3,D^3)$, i.e. there exists a homeomorphism $\psi:M
\lra \R^3$ such that $\psi(X)=D^3$. A { \it $($topological$)$ $3$--ball}
of $M$ is the internal part of a $3$--disk.  We say that a subset $Y$
of $M$ is a {\it $($topological$)$ $2$--disk} if, up to homeomorphism, the
pair $(M,Y)$ is equivalent to $(\R^3,D^2)$. Smooth disks or balls in a
smooth $M$ diffeomorphic to $\R^3$ are defined in the same way by
replacing ``homeomorphism'' with ``diffeomorphism''. Disks and balls in an
arbitrary 3--manifold $W$ are contained, by definition, in some chart
$M$ homeo(diffeo)morphic to $\R^3$.}
\end{defi}

By a {\it domain} $\Omega$ in $\R^3$, we will mean a non--empty connected open
set, which coincides with the interior of its closure in $\R^3$,
i.e. ${\rm Int}\ \overline{\Omega}=\Omega$.  Moreover, throughout the
whole paper, domains will always assumed to be {\it bounded}, whence
with compact closure.

\smallskip

Sometimes it is convenient to identify $\R^3$ with an open subset of
the 3--sphere $S^3=\R^3 \cup \{\infty\}$ via the stereographic
projection from the point ``at infinity''. An open subset
$\Omega\subset S^3$ is a domain if ${\rm Int}\
\overline{\Omega}=\Omega$.  Of course every domain in $S^3$ has
compact closure, and the stereographic projection induces a bijection
between domains in $\R^3$ and domains in $S^3$ whose closure does not
contain the added point $\infty$.

\smallskip

We denote by $\partial \Omega$ the usual (topological) boundary of
$\Omega$, i.e. the set
$$
\partial \Omega=\overline{\Omega}\setminus \Omega \, .
$$
It turns out (see e.g. Remark~\ref{nec-hyp}) that domains with
``wild'' boundary can display pathologi\-cal behaviours that we would
like to exclude from our investigation. We will therefore concentrate
our attention on domains with ``tame'' boundary, carefully
specifying what ``tame'' means in our context.
\smallskip

\subsection{Smooth surfaces.}\label{surfaces}  
We begin by defining the tamest class of domains one could consider.
A {\it smooth surface} $S$ in $\R^3$ is a compact and connected subset
of $\R^3$ such that the following condition holds: for every point $p
\in S$, there exist a neighbourhood $U_p$ of $p$ in $\R^3$ and a
diffeomorphism $\varphi:U_p \lra \R^3$ such that $\varphi (U_p\cap
S)=P$, where $P$ is an affine plane. In other words, $S\subset \R^3$
is a smooth surface if the pair $(\R^3,S)$ is locally modeled, up to
diffeomorphism, on the pair $(\R^3,\R^2)$.  For any system $(x_1,x_2,x_3)$
of linear
coordinates on $\R^3$, for $i=1,2,3$, set $H_i:=\{(x_1,x_2,x_3)\in\R^3 \, | \,
x_i=0\}$. By the Inverse Function Theorem, $S$ is a smooth surface if
and only if it is locally the graph of a real smooth function (defined
on an open subset of some $H_i$).

\begin{prop}\label{divide} 
Every smooth surface $S$ in $\R^3$ disconnects $S^3$ in two domains
$\Omega(S)$ and $\Omega^*(S)$.
\end{prop}

Let us sketch a proof of Proposition~\ref{divide} that uses classical
tools from Differential Topology (exhaustive references for our needs
are, for instance, \cite{MIL} and \cite{HIR}). By the very definition
of surface, if $p$ is a point of $S$, then $S$ disconnects small
neighbourhoods of $p$ into two connected components. Together with the
fact that $S$ is connected, this readily implies that $S^3\setminus S$
consists of at most two connected components. Suppose now, by
contradiction, that $S^3\setminus S$ is connected. Then any closed
interval transverse to $P$ in a local model can be completed in
$S^3\setminus S$ to an embedded smooth circle $f_0:S^1 \lra C_0
\subset S^3$ that transversely intersects $S$ in exactly one
point. Since $S^3$ is simply connected (see Subsection~\ref{iso} for a
brief discussion of such a notion), $f_0$ is smoothly homotopic to an
embedded circle $f_1: S^1 \lra C_1 \subset S^3$ that does not
intersect $S$.  Moreover, we can assume that there exists a smooth
homotopy $F: S^1 \times [0,1] \lra S^3$ between $f_0$ and $f_1$, which
is {transverse} to $S$. Then the set $F^{-1}(S)$ consists of a finite
disjoint union of smooth circles or closed intervals having $F^{-1}(S)
\cap (S^1 \times \{0,1\})$ as set of end--points. In particular,
$F^{-1} (S) \cap (S^1 \times \{0,1\})$ should be given by an {\it
even} number of points, while we know that it consists of just one
point. This gives the desired contradiction.

\smallskip

\noindent {\bf Notation.}  From now on, whenever $S\subset \R^3\subset S^3$
is a smooth surface, we will denote by $\Omega (S)$ and $\Omega^\ast
(S)$ the connected components of $S^3\setminus S$. We will also assume
that $\infty\in \Omega^\ast (S)$, so $\Omega(S)$ is the unique bounded
component of $\R^3\setminus S$, while $\Omega'(S):=
\Omega^*(S)\setminus \{\infty \}$ is the unique unbounded component of
$\R^3\setminus S$. In particular, $\Omega(S)$ is a domain in $\R^3$
and $\partial \Omega (S)= S$.  The {\it local model} of $(\Omega(S),
S)$ at every boundary point is given by $(P_{+},P)$ where $P$ is an
affine hyperplane as above, and $P_{+}\subset \R^3$ is a half--space
bounded by $P$.
\smallskip

\begin{defi}\label{smoothdom}{\rm A domain $\Omega$ in $\R^3$ 
has {\it smooth boundary} if $\partial \Omega$ consists of the
disjoint union of a finite number of smooth surfaces.}
\end{defi}

It readily follows from the definitions that the closure of a domain
with smooth boundary admits a natural structure of compact smooth
manifold with boundary.

The following lemma is an immediate consequence of the previous discussion.

\begin{lem}\label{smoothdom2}
Let $\Omega$ be a domain with smooth boundary. Then we can order the boundary 
surfaces $S_0, S_1,\dots,S_h$ in such a way that: 
\begin{enumerate}
\item
The $\overline{\Omega(S_j)}$'s, $j=1,\ldots,h$, are contained in
$\Omega(S_0)$ and are pairwise disjoint.
\item
$\Omega$ is given by the following intersection:
$$
\Omega = \Omega(S_0) \cap \bigcap_{j=1}^h 
\Omega^*(S_j) \, .
$$
\end{enumerate}
\end{lem}
%\cvd

\smallskip

\subsection{Orientation and tubular neighbourhoods}
Let $S\subset \R^3$ be a smooth surface. We claim that $S$ is {\it orientable}. In
fact, if $\R^3$ is oriented by means of the equivalence class of its
standard basis $(e_1,e_2,e_3)$, then $S$ can be oriented as the
boundary of $\Omega(S)$, via the rule {``first the outgoing normal
vector''}. More explicitly, for each $p \in S$, one can
consistently declare that
a basis $(v_1,v_2)$ of the tangent space $T_p S$ of $S$
at $p$ is positively oriented if and only if
$(n,v_1,v_2)$ is a positively oriented basis
of $\R^3$, where $n$ is a vector orthogonal to $T_p S$ and pointing
outward $\Omega (S)$.

\smallskip

For every $\epsilon >0$, let us define the $\epsilon$--{\it
neighbourhood} $N_\epsilon(S)$ of $S$ in $\R^3$ by setting
$$
N_\epsilon(S):= \{ x\in \R^3 \ | \ {\rm dist}(x,S)\leq \epsilon \} \,
.
$$
If $\epsilon$ is small enough, then the pair $(N_\epsilon(S),S)$
is diffeomorphic to $(S \times [-1,1], S\times \{0\})$. If $r:
N_\epsilon(S)\lra S$ is the natural retraction such that $r(x)$ is the
nearest point to $x$ (such a retraction is well--defined provided that
$\epsilon$ is sufficiently small), then, for every $x \in S$, the set
$r^{-1}(x)$ is a     straight copy of $[-\epsilon,
\epsilon]$. Moreover, $N_\epsilon(S)\cap \overline{\Omega(S)}$ is
mapped onto $S\times [-1,0]$, hence it is     a {\it collar} of $S$ in
$\overline{\Omega(S)}$.  Similarly for     $N_\epsilon(S)\cap
\overline{\Omega'(S)}$. If $C$ is a smoothly embedded     circle in
$\R^3$ and $\epsilon$ is small enough, then $N_\epsilon(C)$ also is a 
tubular neighbourhood of $C$, diffeomorphic to a (closed) solid torus 
$D^2\times S^1$ and having $C$ as {\it core}.

\smallskip

%%%%%%%

\subsection{Link complements.} \label{link}
A {\it link} $L=C_0\cup \ \dots \cup C_h$ in $S^3$ is 
the union of a finite family
of smoothly embedded disjoint circles $C_j$.  If $h=0$, then $L$ is
called a {\it knot}. Suppose that $\infty \in C_0$, hence $A(L)=S^3
\setminus L$ is a connected open set in $\R^3$.  With our definitions,
since $\overline{A(L)}=\R^3$, the internal part of $\overline{A(L)}$
does not coincide with $A(L)$ and $A(L)$ is not a domain. However, to
$L$ there is associated the domain $\compl(L)= S^3 \setminus U(L)$,
where $U(L)$ is the union of small disjoint closed tubular
neighbourhoods of the $C_j$'s. We call $\compl(L)$ {\it
complement--domain of $L$}. The boundary component of $\compl(L)$
corresponding to $C_j$ is a smooth torus $T_j$ and, with the above
notations, $\Omega^*(T_0)$ and $\Omega(T_j)$, $j=1,\dots,h$, are open
solid tori. It is clear that $\compl(L)$ is {\it
homotopically equivalent to $A(L)$} (see e.g. \cite{HATCH} for the definition of homotopy equivalence), hence $\compl(L)$ and $A(L)$ share all the {homotopy type}
invariants (like the fundamental group). A knot $C=C_0$ is {\it
unknotted} if also $\Omega(T_0)$ is a solid torus or, equivalently, if
$C$ bounds a 2--disk of $S^3$. A link has {\it geometrically
unlinked components} if its components are contained in pairwise
disjoint 3--disks of $S^3$. A link is {\it trivial} if it has geometrically unlinked unknotted components.

Suppose now that $\infty \not\in L$, i.e. consider $L$ as a link of
$\R^3$. We use the symbol $U(L)$ again to indicate the union of small
disjoint closed tubular neighbourhoods of the $C_j$'s in
$\R^3$. Choose a smooth $3$--ball $B$ of $\R^3$ containing $U(L)$ and
define $\B(L):=B \setminus U(L)$. We call $\B(L)$ {\it box--domain}
of~$L$. Any rigid motion of $S^3$ that takes $L$ onto
a link $L'$ containing the point at infinity establishes a 
diffeomorphism between the box--domain 
$\B(L)$ and the complement--domain $\compl(L')$ with a $3$--disk
removed.

The reader observes that the complement-- and the box--domains of a
link are well--defined, up to diffeomorphism (up to ambient isotopy
indeed).

%%%%%%%

\subsection{Cutting along surfaces.}                                   
Let $\Omega$ be a domain with smooth boundary. A {\it properly
embedded surface $\Sigma$ in $(\overline {\Omega},\partial \Omega)$}
is a {compact and connected} subset of $\overline{\Omega}$ such that:

\begin{enumerate}
\item On $\Sigma \setminus \partial \Omega$, $\Sigma$ has the same
local model of a smooth surface. 

\item If $\Sigma \cap \partial \Omega \neq \emptyset$, then at every
point of this intersection, up to local diffeomorphism, the triple
$(\overline {\Omega},\partial \Omega, \Sigma)$ is equivalent to the
{local model} $(P_+,P,T_+)$, where $(P_+,P)$ are as in
Subsection~\ref{surfaces}, and $T_+=T\cap P_+$, $T$ being a plane
orthogonal to $P$.  It follows that $\Sigma$ is a {\it smooth surface
with boundary} $\partial \Sigma = \Sigma \cap \partial \Omega$. This
boundary is a (not necessarily connected) smooth curve embedded in
$\partial \Omega$.

\item $(\Sigma, \partial \Sigma)$ admits a {\it bicollar} in
$(\overline {\Omega},\partial \Omega)$, i.e.~there exists a closed
neighbourhood $U$ of $\Sigma$ in $\overline\Omega$ such that
$(U,U\cap\partial\Omega)$ is diffeomorphic to $(\Sigma \times [-1,1],
(\partial \Sigma)\times [-1,1])$, via a diffeomorphism sending each
point $x \in \Sigma$ into $(x,0) \in \Sigma\times \{0\}$.  It is not
hard to see that the existence of a bicollar is equivalent to the fact
that $\Sigma$ is orientable.  Any orientation on $\Sigma$ induces an
orientation on $\partial \Sigma$, via the rule ``first the outgoing
normal vector'' mentioned above.
\end{enumerate}

Let $\Sigma$ be properly embedded in $(\overline {\Omega},\partial
\Omega)$.
% and assume that $\partial \Sigma$ is non--empty.  
Then the
result $\Omega_C(\Sigma)$ of the {\it cut/open} operation along
$\Sigma$ consists in taking the internal part in $\R^3$ of the
complement in $\overline {\Omega}$ of a bicollar of $(\Sigma,\partial
\Sigma)$. In general, $\Omega_C(\Sigma)$ is not connected. However,
every connected component of $\Omega_C(\Sigma)$ is a domain.  The
boundary of $\Omega_C (\Sigma)$ is no longer smooth, because some
{corner lines} arise along $\partial \Sigma$. However, by means of a
standard ``rounding the corners'' procedure, we can assume that {the
class of domains with smooth boundary is closed under the cut/open
operation}.

\begin{remark}\label{on-cut}{\rm
 A more direct way to cut should be by taking $A(\Sigma)= \Omega
\setminus \Sigma$. The components of $A(\Sigma)$ are not domains in
general. On the other hand, each component of $\Omega_C(\Sigma)$ is
contained in and is homotopically equivalent to one component of
$A(\Sigma)$.  This establishes a bijection between these two sets of
components, and corresponding components of $A(\Sigma)$ and
$\Omega_C(\Sigma)$ share all the homotopy type invariants.}
\end{remark}

%\smallskip

\begin{exa}\label{Seifert}
{\rm Given a knot $K$ in $S^3$, a {\it Seifert surface} of $K$ is a
connected orientable smoothly embedded surface $S$ with boundary equal
to $K$. Every knot has a Seifert surface (see \cite{rolfsen}). Given
the domain $\compl(K)$ as in Subsection \ref{link}, we can assume that
such a surface $S$ is transverse to the boundary torus along a {\it
preferred longitude} parallel to $K$ (it is well-known that the
isotopy class of this preferred longitude does not depend on the
chosen Seifert surface -- see Remark~\ref{uniqueslope}).
Hence,
$\Sigma:= S \cap \compl(K)$ is properly embedded in $\compl(K)$ and
the corresponding cut/open domain $(\compl(K))_C(\Sigma)$, being
connected, is a domain.}
\end{exa} 

%%%%%%%

\subsection{Locally flat boundary} 
In order to perform constructions and develop arguments which use
tools from Differential Topology, it is very convenient to work with
smooth boundaries. Such a choice allows us, for instance, to exploit
the powerful notion of {\it transversality}. We have already used such
a notion in the proof of Proposition~\ref{divide} sketched
above. Moreover, using transversality, we will be able to approach in
an elementary, geometric and quite ``primitive'' way some fundamental
results about {\it duality} (such results are usually esta\-blished in
more general settings by using more sophisticated tools from Algebraic
Topology). On the other hand, people dealing with Helmholtz domains
usually work with boundaries of weaker classes of regularity, in
particular with boundary that are local graphs of Lipschitz
functions. In this case, the domain is said to have {\it Lipschitz
boundary}. A natural way to deal with more general topological
boundaries, keeping nevertheless the same qualitative local pictures,
consists in considering triples $(\overline {\Omega}, \partial \Omega,
\Sigma)$ that admit everywhere the same (suitable) local models of the
smooth case, providing that we replace ``up to local diffeomorphism''
with {``up to local homeomorphism''}. Such {topological} triples are
called {\it locally flat}. Note that, according to these definitions,
our topological disks in 3--manifolds are locally flat. The following
lemma is immediate.

\begin{lem}\label{loc-graph} 
A compact connected subset of $\R^3$, which is locally the graph
of {\rm continuous} fun\-ctions, is a locally flat surface.
\end{lem}
%\cvd
\smallskip

Several deep fundamental results of 3--dimensional Geometric Topology
\cite{MOI, bing, brown2} imply that, up to homeomorphism, there is not
a real difference between the smooth and the locally flat topological
case:

\begin{teo}\label{smoothing} 
For every locally flat triple $(\overline {\Omega},\partial \Omega,
\Sigma)$, the following statements hold.
\begin{itemize}
\item[$(1)$] {\bf Triangulation.} There is a homeomorphism $t:\R^3
\lra \R^3$ that maps the given triple onto a {\rm polyhedral} triple
$($i.e.  the piecewise linear realization in $\R^3$ of a finite {\rm
simplicial complex} with distinguished subcomplexes$)$.
\item[$(2)$] {\bf Smoothing.} There is a homeomorphism $s:\R^3 \lra \R^3$ that
maps the given triple onto a smooth one.
\end{itemize}
\end{teo}
%\cvd

\smallskip

Summarizing:

\smallskip

{\it In order to study the geometric topology of {\rm arbitrary}
locally flat topological triples, it is not restrictive to consider
only smooth ones. Moreover, if useful, we can use also tools from
$3$--dimensional {\rm Polyhedral (PL) Geometry.}}

\smallskip

\subsection{Isotopy, homotopy and homology}\label{iso}
Before entering the main part of our work, we would like to give a
brief and intuitive description of some concepts that will be
extensively used throughout the paper (they will be treated a bit more
formally in Sections \ref{algtop1:subsec} and \ref{algtop2:subsec}).
Let $M$ be a smooth connected $n$--manifold with (possibly empty)
boundary (for our purposes, it is sufficient to consider the cases in
which $M$ is a 3--dimensional domain as above or the whole spaces
$\R^3,S^3$, or a smooth surface).  Two smooth simple oriented loops
$C_0,C_1 \subset M$ are \emph{isotopic} if they are related by a
smooth isotopy, i.e.~by a smooth map $F\colon S^1 \times [0,1]\lra M$
such that, if $F_t:= F(\cdot , t)\colon S^1\lra M$, then $F_0,F_1$ are
oriented parameterizations of $C_0,C_1$ respectively, and $F_t$ is a
smooth embedding for every $t\in [0,1]$.  In other words, $C_0$ is
isotopic to $C_1$ if it can be smoothly deformed into $C_1$ without
crossing itself.

A \emph{homotopy} between $C_0$ and $C_1$ is just the same as an
isotopy, provided that we do \emph{not} require $F_t$ to be an
embedding for every $t$. More precisely, if $C_0,C_1$ are
\emph{continuous} (possibly non--injective) loops of $M$, we say that
$C_0$ is homotopic to $C_1$ if it can be taken into $C_1$ by a
\emph{continuous} deformation along which non--injectivity phenomena
such as self--crossings are allowed.  In particular, $C_0$ is
homotopically trivial if it is homotopic to a constant loop, or,
equivalently, if a parametrization of $C_0$ can be extended to a
continuous map from the $2$--disk $D^2$ to $M$ (where we are
identifying $S^1$ with $\partial D^2$). The manifold $M$ is {\it
simply connected} if (it is connected and) every loop in $M$ is
homotopically trivial.  It is well--known (and very easy) that $\R^3$
and $S^3$ are simply connected, while by the very definition
non--trivial knots in $S^3$ provide examples of loops that are not
isotopic to the unknot. Recall that unknotted knots can be
characterized as those knots which bound a $2$--disk. 

More in general,
let us define a \emph{$1$--cycle} (\emph{with integer coefficients})
in $M$ as the union $L$ of a finite number of (not necessarily
embedded nor disjoint) oriented loops in $M$. We say that $L$ is a
\emph{boundary} if there exist an oriented (possibly disconnected)
surface with boundary $S$ and a continuous map $f\colon S\lra M$ such
that the restriction of $f$ to the boundary of $S$ defines an
orientation--preserving parameterization of $L$ (the orientation of
$S$ canonically induces an orientation of $\partial S$ also in the
topological setting): with a slight abuse, in this case, we say that
$L$ \emph{bounds} $f(S)$.  Of course, knots and links in $S^3$ are
particular instances of $1$--cycles in $S^3$, and every knot is a
boundary, since it bounds a (possibly singular) $2$--disk, or a
Seifert surface. If $L,L'$ are $1$--cycles in $M$ and $-L'$ is the
$1$--cycle obtained by reversing all the orientations of the loops of
$L'$, we say that $L$ is \emph{homologous} to $L'$ if the $1$--cycle
$L\cup -L'$ is a boundary, and that $L$ is homologically trivial if it
bounds or, equivalently, if it is homologous to the empty
$1$--cycle. It readily follows from the definitions that homotopic
loops define homologous $1$--cycles. The space of equivalence classes
of $1$--cycles is called \emph{singular $1$--homology module of $M$}
(\emph{with integer coefficients}) and it is usually denoted by $H_1
(M;\Z)$. The union of cycles induces a sum on $H_1 (M;\Z)$, which is
therefore an Abelian group. It is not difficult to show that, since
$M$ is connected, every $1$--cycle in $M$ is homologous to a single
loop, and this readily implies that, if $M$ is simply connected, then
$H_1 (M;\Z)=0$. The converse statement is not true in general (see
Remark~\ref{homologysphere}), but turns out to hold for tame domains
in $\R^3$ (see~Corollary~\ref{cor-simple}). Note, however, that even
if $M=\Omega$ is a domain in $S^3$ with locally flat boundary, then
there may exists a loop of $M$ which is homologically trivial, but not
homotopically trivial: if $K \subset S^3$ is a non--trivial knot with
complement--domain $\compl(K)$, then a Seifert surface $\Sigma$ for
$K$ defines a preferred longitude $\gamma=\Sigma \cap \partial
\compl(K) \subset \partial \compl(K)$. Such a longitude bounds the
surface with boudary $\Sigma \cap \overline{\compl(K)}$ and is
therefore homologically trivial in $\overline{\compl(K)}$. However, as
a consequence of the classical Dehn's Lemma (see \cite[p. 101]{rolfsen}),
if $\gamma$ were homotopically trivial in $\overline{\compl(K)}$, it
would bound a (embedded locally flat) $2$--disk in
$\overline{\compl(K)}$, and this would imply in turn that $K$ is
trivial, a contradiction.  

The singular $2$--homology module of $M$
can be described in a similar way as the set of equivalence classes  
of maps of compact smooth oriented
(possibly disconnected) surfaces in $M$, up
to 3--dimensional ``bordism''.  A nice, non--trivial fact in
the situations of our interest, is that every 1-- or 2--homology class
can be represented by {\it submanifolds} (i.e the above maps are
embeddings), and that also the bordisms between homologically
equivalent submanifolds can be realized by submanifolds.  In the
polyhedral setting, this is a consequence of {\it Kneser's method}
(1924) for eliminating singularities (see \cite[p. 32]{Marin}). By
Theorem \ref{smoothing} (or even by classical results within the 
smooth framework), this
holds also in the smooth case.

%%%%%%%%%%%%%%%%%%%%%

\section{Simple domains}\label{simple}

Let $\Omega$ be a domain. In theoretic and applied literature
about Helmholtz domains, two main notions are employed in
order to specify the way $\Omega$ is {``simple''}:
\smallskip

(a) {\it $\Omega$ is simply connected} (i.e.~has trivial fundamental group).
\smallskip

(b) {\it Every curl--free smooth vector field on $\Omega$ is the gradient of
a smooth function on $\Omega$.}

\smallskip

\noindent
Other related conditions will be considered in Corollary~\ref{cor-simple}.

\smallskip

It is widely well--known (see anyway the corollary just mentioned) that
$$
{\rm (a)} \Longrightarrow {\rm (b)} \, .
$$ We are going to discuss presently the converse implication, which
seems to have risen some misunderstandings (see Example~\ref{Vou-Bin}
below).

\smallskip

\subsection{Vector fields, differential forms and de Rham cohomology.}
\label{duality:subsec} 
We begin by reformulating condition (b) more conveniently in terms of
differential forms. It is well--known from Linear Algebra that every
non--degenerate scalar product $\langle \ , \, \rangle$ on a finite
dimensional real vector space $V$ determines an isomorphism $\psi: V
\lra V^*$ between $V$ and its {\it dual space} $V^*:= {\rm
Hom}(V,\R)$, by the formula $\psi(v)(w)=\langle v,w \rangle$, for
every $v,w \in V$. A Riemannian metric on a smooth manifold $M$ is
just a smooth field $\{\langle \ , \, \rangle_p\}_{p\in M}$ of
positive definite (hence non--degenerate) scalar products on the
tangent spaces $T_pM$. The same formula applied pointwise at every
point $p$ of $M$ determines a canonical isomorphism between the space
of smooth tangent vector fields and the space of smooth 1--forms on
$M$ (from now on, even when not explicitly stated,
differential forms will always assumed to be smooth). Let us apply this general fact to the standard flat Riemannian
metric $ds^2=dx_1^2+dx_2^2+dx_3^2$ on $\R^3$ (and to its restriction
to any domain). In practice, if $V=(V_1,V_2,V_3)$ is a smooth vector
field on a domain $\Omega$, then $\omega := \sum_{j=1}^{3} V_jdx_j$ is
the associated 1--form. The differential of $\omega$ is the 2--form
$$
d\omega=\left(-\frac{\partial V_2}{\partial x_3}+\frac{\partial V_3}
{\partial x_2}\right) \, dx_2 \wedge dx_3
-
\left(\frac{\partial V_1}{\partial x_3}-\frac{\partial V_3}
{\partial x_1}\right) \, dx_1 \wedge dx_3
+
\left(-\frac{\partial V_1}{\partial x_2}+\frac{\partial V_2}
{\partial x_1}\right) \, dx_1 \wedge dx_2 \, .
$$ 
Since 
$$
{\rm curl}(V)=\left(-\frac{\partial V_2}{\partial x_3}+\frac{\partial V_3}
{\partial x_2}, \, \frac{\partial V_1}{\partial x_3}-\frac{\partial V_3}
{\partial x_1}, \, -\frac{\partial V_1}{\partial x_2}+\frac{\partial V_2}
{\partial x_1}\right) \, ,
$$ 
$V$ is curl--free if and only if $d\omega=0$.

If $f\colon \Omega\lra \R$ is a smooth function, 
the differential of $f$ is the 1--form 
$$
df=\sum_{j=1}^3\frac{\partial f}{\partial x_j} \, dx_j \, .
$$ By the very definitions, the gradient $\nabla f$ corresponds to
$df$, via the above canonical isomorphism determined by $ds^2$.

A 1--form is \emph{closed} if its differential vanishes, and it is
\emph{exact} if it is the differential of a smooth function. Since
$d(df)=0$ for every smooth function $f$ (or, equivalently, every
gradient field is curl--free), every exact 1--form is closed. If
$\Omega$ is a domain, then the first de Rham cohomology group
$H^1_{DR} (\Omega)$ is defined as the quotient vector space of closed
1--forms defined on $\Omega$ modulo exact 1--forms defined on
$\Omega$. Condition (b) above is then equivalent to condition

\smallskip

($\rm b'$) {\it Every closed $1$--form on $\Omega$ is exact,
i.e. $H^1_{\rm DR}(\Omega)=0$.}

\smallskip

This already shows that condition $\rm (b')$ only depends on the
differential structure of $\Omega$, and it is not necessary to drag
the Riemannian metric in, like one actually does in (b). Moreover, as
a very particular case of de Rham's Theorem (see e.g. \cite{BOT}), we
know that
$$
H^1 _{\rm DR}(\Omega) \cong H^1(\Omega;\R) \, ,
$$ where the vector space on the right--hand side is the {\it singular
$1$--cohomology module with real coefficients}, which is a topological
(homotopic type indeed) invariant. Hence, we have a new reformulation
of (b) in terms of basic notions taken from Algebraic Topology (an
exhaustive reference for our needs is \cite{HATCH}):

\smallskip

$\rm (b'')$  $H^1(\Omega;\R)=0$.

\smallskip

We are now ready to state the main result of this section, which
provides an easy chara\-cterization of simple domains in $\R^3$. 
We keep notations from Lemma~\ref{smoothdom2} and defer the proof to
Subsection~\ref{mainproof:subsec}.

\begin{teo}\label{simple-top} 
Let $\Omega$ be a domain of $\R^3$ with locally flat boundary such
that $H^1(\Omega;\R)=0$. Then, for every $j \in \{0,1,\dots,h\}$, both
$\Omega(S_j)$ and $\Omega^*(S_j)$ are $3$--balls of $S^3$ bounded by
the locally flat $2$--sphere $S_j$. In particular, $\Omega$ is simply
connected.
\end{teo}

\smallskip

Such a result can be rephrased as follows:

\smallskip

{\it Every domain of $\R^3$ with locally flat boundary and with
$H^1(\Omega;\R)=0$ consists of an {``external''} $3$--ball with some
$($a finite number indeed$)$ {``internal''} pairwise disjoint
$3$--disks removed}.

\smallskip

Singular homology and singular cohomology with real and integer
coefficients are closely related to each other by the Universal
Coefficient Theorem (see~ e.g.~\cite{HATCH}).  We now list two
easy consequences of this classical result, which will prove useful for
establishing the equivalence between the different definitions of
simple domain described in the following corollary.  More details can
be found in Subsections~\ref{algtop1:subsec} and~\ref{algtop2:subsec}.

Let $X$ be any topological space. Denote by $H_1 (X;\R)$ the singular
1--homology module of $X$ with real coefficients, and 
recall that $H_1
(X;\Z)$ is the singular 1--homology module of $X$ with integer
coefficients. Then the Universal Coefficient Theorem provides the
following canonical isomorphisms
$$
H^1(X;\R)\cong {\rm Hom}(H_1(X;\R),\R),
\qquad 
H_1 (X;\R)\cong H_1 (X;\Z)\otimes \R \, .
$$ 

\begin{cor}\label{cor-simple} 
Let $\Omega$ be a domain with locally flat boundary. Then the
following properties are equivalent:
\begin{itemize}
\item[(a)] $\Omega$ is simply connected.
\item[(b)] Every curl--free smooth vector field on $\Omega$ is the gradient of
a smooth function.
\item[($\rm b''$)] $H^1(\Omega;\R)=0$.
\item[(c)]  $H_1(\Omega;\Z)=0$.
\item[(d)] $H_1 (\Omega;\R)=0$.
\item[(e)] For every curl--free smooth vector field $V$ and every
divergence--free smooth vector field $W$ on $\Omega$ with compact
support, the integral $\int_{\Omega}V \bullet W \, dx$ is null, where
$V \bullet W:=\sum_{j=1}^3 V_j \cdot W_j$ if $V=(V_1,V_2,V_3)$ and
$W=(W_1,W_2,W_3)$.
\end{itemize}

Moreover, if $\Omega$ has Lipschitz boundary, then we
can add the following equivalent condition to the list:
\begin{itemize}
\item[(f)] Every vector field in $(L^2(\Omega))^3$ with null
distributional curl is the weak gradient of a function in the Sobolev
space $H^1(\Omega)$ $($here $H^1(\Omega)$ denotes the set of all
elements of $L^2(\Omega)$ having weak gradient in $(L^2(\Omega))^3)$.
\end{itemize}
\end{cor}
\Dim As observed in Subsection~\ref{iso}, if $\Omega$ is simply
connected, then every $1$--cycle in $\Omega$ is a boundary, so $H_1
(\Omega;\Z)=0$. As a consequence of the Universal Coefficient Theorem,
we have then $H_1 (\Omega;\R)=0$ and $H^1 (\Omega;\R)=0$. We have thus
proved that
$$
{\rm (a)} \Longrightarrow {\rm (c)} \Longrightarrow {\rm (d)}
\Longrightarrow {\rm (b'')} \ 
\left(\Longleftrightarrow {\rm (b)}\right)
\, .
$$ On the other hand, Theorem~\ref{simple-top} ensures that $\rm
(b'')$ implies (a). We have thus proved that the first five conditions
are equivalent to each other.

If (b) holds, then (e) follows immediately from the Green
formula. Suppose now that (e) holds, let $V$ be a curl--free smooth
vector field on $\Omega$ and let $\omega$ be the $1$--form
corresponding to $V$ via the duality described above. 
Let now $\varphi$ be any fixed compactly supported closed
$2$--form on $\Omega$.  As a direct consequence of Stokes' Theorem,
the map which associates to every class $[\psi]\in H^1_{DR} (\Omega)$
the real number
$$
\int_\Omega \psi\wedge \varphi
$$ is well--defined and determines therefore a linear map
$f_\varphi\colon H^1_{DR} (\Omega) \lra \R$.  Now a classical result
in de Rham Cohomology Theory (see e.g.~\cite[p. 44]{BOT}) ensures that
every linear map $H^1_{DR} (\Omega) \lra \R$ arises in this way,
i.e. it is of the form $f_{\varphi}$ for some closed compactly
supported 2--form $\varphi$.  Therefore condition (e) translates into
the fact that every linear map $H^1_{DR} (\Omega)\lra \R$ vanishes on
the cohomology class $[\omega]$ of $\omega$, and this readily implies
that $[\omega]=0$, i.e. $\omega$ is exact. This is in turn equivalent
to the fact that $V$ is the gradient of a smooth function.

Finally, $\rm (b'') \Longleftrightarrow (f)$  is immediate from the version
of de Rham's Theorem given in \cite[Assertion (11.7), p. 85]{mitrea}.  \cvd

\subsection{A fallacious counterexample}
Before going into the proof of Theorem~\ref{simple-top}, we discuss a
{\it fake counterexample} to (b)$\Longrightarrow$(a).

\begin{exa}\label{Vou-Bin}{\rm 
This is a fallacious example given by A. Vourdas and K. J. Binns in
their response to R. Kotiuga in the correspondence \cite[p. 232]{KVB}
(see also \cite{BVB}, \cite[Subsection~2.1]{KKB} and \cite[Section
1]{KKB2}). We refer to Remark \ref{on-cut}, Subsections \ref{link} and
\ref{Seifert}.  Let $C$ be the oriented trefoil knot of $\R^3$ and let $\Sigma$
be the Seifert surface of $C$ drawn in Figure \ref{fig:unique} (on the
left). Denote by $\Omega_C(\Sigma)$ the domain of $\R^3$ obtained
by applying to the complement--domain $\compl(C)$ of $C$ the cut/open
operation along $\Sigma$.

%\vspace{-1em}

\begin{figure}[htbp]
\begin{center}
 \includegraphics[height=5.52cm]{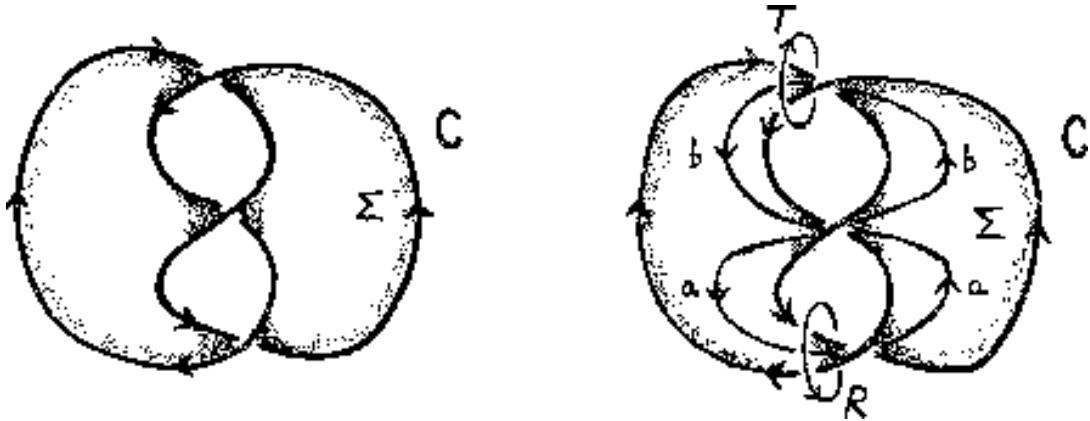}
\caption{\label{fig:unique} The trefoil knot with one of its Seifert surfaces.}
\end{center}
\end{figure}

In \cite[p. 232]{KVB}, the authors assert that $H_1(S^3\setminus
\Sigma;\R)=0$ (equivalently $H_1(\Omega_C(\Sigma);\R)=0$), but that
$S^3\setminus \Sigma$ (equivalently, $\Omega_C(\Sigma)$) is not simply
connected.  The first claim is wrong.  In fact, consider the two
oriented loops $a$ and $b$ contained in $\Sigma$ and the two oriented
loops $R$ and $T$ contained in $\Omega_C(\Sigma)$ drawn in
Figure~\ref{fig:unique} (on the right). The surface $\Sigma$ is
homeomorphic to a torus minus an open $2$--ball, and the homology
classes of $a$ and of $b$ in $\Sigma$ form a basis of $H_1(\Sigma;\R)$
(see also Figure 8.12 of \cite[p. 243]{bossavit} to visualize these
facts). Lefschetz Duali\-ty Theorem immediately implies that the
homology classes of $R$ and of $T$ form a basis of
$H_1(\Omega_C(\Sigma);\R)$. In particular, this last space is
non--trivial. Moreover, the trefoil knot is an example of {\it fibred
knot} having the given Seifert surface as a fibre (this is carefully
described in \cite[p. 327]{rolfsen})).  Hence, $\Omega_C(\Sigma)$ is
homeomorphic to $(\Sigma \setminus \partial \Sigma)\times (0,1)$ and
has therefore the same homotopy type of $\Sigma$. Note that this fact
confirms the above claim that $H_1(\Omega_C(\Sigma);\R)$ and
$H_1(\Sigma;\R)$ are isomorphic.

The first argument above can be rephrased in a more physical
fashion. Suppose $a$ is an ideally thin conductor, carrying a current
of unitary intensity. Let $\HH_a$ be the corresponding magnetic
field. The restriction $\HH^{\prime}_a$ of $\HH_a$ to $S^3\setminus
\Sigma$ is a curl--free smooth vector field, which does not have any
scalar potential. In fact, the circulation of $\HH^{\prime}_a$ along
$R$ is $1$. In particular, by Stokes' Theorem, the homology class of
$R$ in $S^3\setminus \Sigma$ is not null. Similar considerations can be
repeated for $b$ and $T$. 

We believe that the following observation 
contains a possible source of
this mistake. In Figure~\ref{fig:boundary}, it is drawn a compact
connected orientable surface $B$ of $\R^3$ with boundary $R$ contained
in $S^3 \setminus C$ (see also Figure 8.13 of
\cite[p. 244]{bossavit}).  The existence of such a surface implies
that $R$ represents the null homology class in $H_1(S^3 \setminus C;
\R)$. Then the restriction to $S^3\setminus \Sigma$ of any curl--free
smooth vector field defined on the {\it whole} of $S^3\setminus C$ has
null circulation along $R$. On the other hand, not every curl--free
smooth vector fields on $S^3\setminus \Sigma$ can be extended to
$S^3\setminus C$. Note also that the surface $B$ intersects in an
essential way the Seifert surface $\Sigma$. These facts explain why
the homology class of $R$ in $S^3 \setminus C$ is null, while the
homology class of $R$ in $S^3\setminus \Sigma$ is not. We will
elaborate this remark in Section~\ref{weakly-Helmholtz} below.}
\end{exa}  

%\vspace{-1em}

\begin{figure}[htbp]
\begin{center}
 \includegraphics[height=5.1cm]{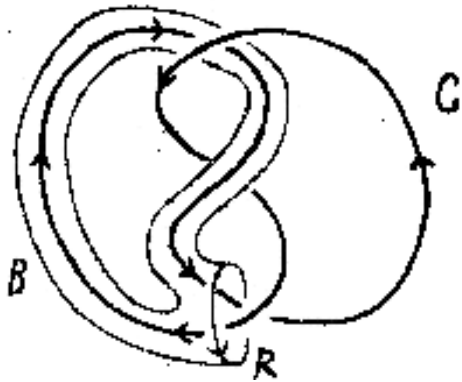}
\caption{\label{fig:boundary} A null homologous cycle.}
\end{center}
\end{figure}

\begin{figure}[htbp]
\begin{center}
 \includegraphics[height=10cm]{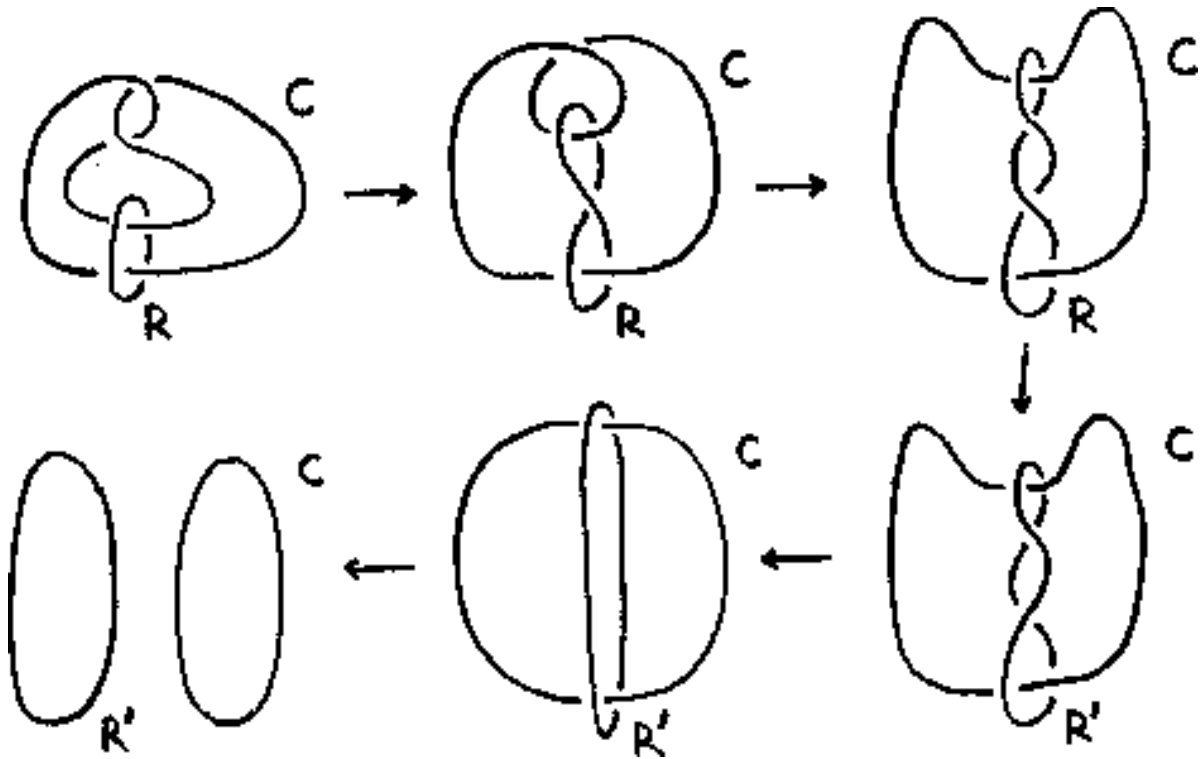}
\vspace{-1em}
\caption{\label{fig:white} Homotoping $R$ to a 
trivial knot in $\R^3\setminus C$.}
\end{center}
\end{figure}

\begin{exa}\label{nuovoesempio}
{\rm In their discussion about the relationship between homotopy and
homology, Vourdas and Binns also consider the case of the
\emph{Whitehead link} (see Figure~\ref{fig:white} above on the left).
With notations as in Figure~\ref{fig:white}, they claim that the loop
$R$ is homologically trivial and homotopically non--trivial in the
complement of $C$ (see \cite{BVB}). On the contrary, the sequence of
moves described in Figure~\ref{fig:white} shows that $R$ is homotopic
(in the complement of $C$) to a loop $R'$ which is clearly
null--homotopic. In fact, as discussed in Subsection~\ref{iso}, since
$R,R'$ are loops in $\R^3\setminus C$ and $R$ can be continuosly
deformed into $R'$ without crossing $C$ (but crossing itself!), then
$R$ and $R'$ are homotopic in $\R^3 \setminus C$. This implies, in
particular, that $R$ bounds a singular $2$--disk in $\R^3\setminus
C$. In fact, since $R$ and $C$ are not geometrically unlinked, $R$
\emph{cannot} bound an \emph{embedded} locally flat $2$--disk in
$\R^3\setminus C$.  As a consequence, it can be shown that $R$ and
$R'$ are not isotopic in $\R^3\setminus C$.}
\end{exa}

\subsection{Elementary results about the algebraic topology of domains}
\label{algtop1:subsec}
Let $M$ be a compact smooth manifold. We say that $M$ is {\it closed} if its
boundary is empty. By the classical {\it Morse theory} (see \cite{MIL2},
\cite{HIR}), if $M$ is closed, then it has the homotopy type of a finite
CW complex of dimension $m = {\rm dim}\ M$, which can be constructed by
means of any Morse function on $M$.  If $M$ is connected with {\it
non--empty} boundary, then it has the homotopy type of a CW complex of
dimension $<m$. This can be realized by means of any Morse function
$f:(M,\partial M)\lra ([0,1],\{1\})$ without local maxima. The same
facts hold if $M$ is polyhedral. One can get a unified treatment by
reformulating Morse theory in terms of {\it handle decompositions theory},
which makes sense also in the polyhedral setting (see \cite{MIL3},
\cite{R-S}). By Theorem~\ref{smoothing}, in our favourite case of
spatial domains, we can adopt both points of view.

Since $\R$ is a field, an easy application of the Universal
Coefficient Theorem for cohomology shows that, for every $k \in
\mathbb{N}$, the singular $k$--cohomology module $H^k(M;\R)$ of $M$
with coefficients in $\R$ is isomorphic to the dual space ${\rm
Hom}(H_k(M;\R),\R)$ of the corresponding singular homology module
$H_k(M;\R)$. Moreover, compactness of $M$ implies that, for every $k
\in \mathbb{N}$, the \emph{$k$--th Betti number} $b_k (M):= \dim
H_k(M;\R)$ of $M$ is finite, whence equal to $\dim H^k(M;\R)$. In
fact, by using the fundamental isomorphism between {\it cellular} (or
{\it simplicial}) and singular homologies, it follows that $\dim
H_k(M;\R)$ is finite and vanishes for every $k> \dim M$. Similar
results also hold for homology and cohomology with integer
coefficients: $H_n (M;\Z)$ and $H^n (M;\Z)$ are finitely generated for
every $n\in\mathbb{N}$ and trivial for $n>\dim M$. Hence, if we denote
by $T_n (M)$ the submodule of finite--order elements of $H_n (M,\Z)$,
then $T_n (M)$ is finite and
$$
H_n(M;\Z)= \big(H_n(M;\Z)/T_n(M)\big) \oplus T_n(M) \, .
$$ Being finitely generated and torsion--free, the quotient $H_n
(M;\Z)/T_n(M)$ is isomorphic to $\Z^r$ for some $r \geq 0$; such a $r$
will be called the \emph{rank} of $H_n (M;\Z)$ and will be denoted by
$r_n (M)$.

Since $\R$ is a field, the Universal Coefficient Theorem for homology
ensures that $H_n(M;\R)=H_n(M;\Z)\otimes \R$, and this implies in turn
that $r_n (M)=b_n (M)$. Let us now recall the definition of the {\it
Euler--Poincar\'e characteristic} $\chi(M)$ of $M$:
$$
\chi(M):= \sum_{n=0}^{\dim  M} (-1)^n b_n(M) \, .
$$ It is well--known that, if $c_n$ is the number of $n$--cells
($n$--simplexes) of any finite CW complex homotopy equivalent to (any
triangulation of) $M$, then $\chi (M)$ admits the following
combinatorial description:
$$
\chi(M)= \sum_{n=0}^{\dim  M} (-1)^n c_n \, .
$$

\smallskip

We now list some elementary results that will prove useful later.

\smallskip 

(1) Assume that $M$ is connected. Then $b_0(M)=1$. If $\dim M=m$ and $M$
has {\it non--empty} boundary, then $b_m(M)=0$.  The last claim follows
from the above--mentioned fact that $M$ has the homotopy type of a CW
complex of stricly smaller dimension.

\smallskip 

(2) If $M$ is a closed manifold of {\it odd} dimension $m=2n+1$, then
$\chi(M)=0$. In fact, by using the ``dual'' CW complexes associated to
$f$ and $-f$, where $f$ is a suitable Morse function on $M$, one
realizes that the respective numbers of cells verify the relations
$c_i = c^*_{m-i}$, and hence the result easily follows from the combinatorial
formula for the characteristic. If $M$ is triangulated, one can use the
{\it dual cell decomposition} of a given triangulation. This is a primitive
manifestation of the {\it Poincar\'e duality} on $M$.

\smallskip 

(3) If $M$ is a connected manifold with {\it non--empty} boundary
$\partial M$, then we can construct the {\it double} $D(M)$ of $M$, by
glueing two copies of $M$ along their boundaries via the identity
map. Then $D(M)$ is closed and
$$
\chi (D(M)) = 2\chi(M) - \chi (\partial M) \, .
$$ In the case of triangulable manifolds (like spatial domains), the
latter equality follows easily by considering a triangulation of
$(M,\partial M)$, that induces a triangulation of the double, and by
using the combinatorial formula for $\chi$.  Hence, if $\dim M$ is
odd, then $\chi (\partial M) = 2\chi (M)$ is even. Moreover, we
observe that
$$
\chi (\partial M) = \sum_i \chi (S_i) \, ,
$$
where the $S_i$'s are the boundary components of $M$.

\smallskip

Let us now specialize to domains.

\smallskip

(4) As already mentioned, if $\Omega \subset \R^3$ is a domain with
smooth boundary, then $\Omega$ is homotopically equivalent to
$\overline{\Omega}$, so $b_n (\Omega)=b_n (\overline{\Omega})$ for
every $n \in \mathbb{N}$. Since $\overline{\Omega}$ is a compact smooth
$3$-manifold with non--empty boundary, we deduce from point~(1) above that 
$$
\chi (\Omega)= \chi (\overline{\Omega}) = 
1- b_1 (\Omega) + b_2 (\Omega) \, .
$$

\smallskip

(5) If $M=S$ is a smooth surface in $\R^3$, then $b_0(S)=1=b_2(S)$, and
    $S$ bounds $\overline{\Omega(S)}$. In particular, by point
    (3) above, $b_1(S)=2-\chi(S)$ is even. The non--negative integer
$$
g(S):=\frac{b_1(S)}{2}
$$ is called {\it genus} of $S$. A basic classification theorem of
orientable surfaces (see \cite{HIR}) says that two compact orientable
surfaces are diffeomorphic if and only if they have the same genus. In
particular, $S$ is a smooth 2--sphere if and only if $g(S)=0$.

\smallskip

(6) If $\Omega\subset \R^3$ is a domain whose boundary consists of the
    disjoint union of smooth surfaces $S_0,\ldots,S_h$, then, by
    points (3) and (5) above, it holds:
$$ \chi(\Omega)=\frac{\chi(\partial
M)}{2}=\frac{1}{2}\sum_{i=0}^h\chi(S_i)=\frac{1}{2}\sum_{i=0}^h \big(
2-2g(S_i) \big)=h+1 - \sum_{i=0}^h g(S_i) \, .
$$

%%%

\subsection{Proof of Theorem \ref{simple-top}}\label{mainproof:subsec}
Let $\Omega$ be a domain with locally flat boundary such that
$H^1 (\Omega; \R)=0$.  We know that it is not restrictive to
assume that $\Omega$ has smooth boundary. We denote by $S_0,\ldots,
S_h$ the boundary components of $\partial \Omega$, keeping notations
from Lemma \ref{smoothdom2}.

Let us set $b_1:=b_1 (\Omega)$, $b_2:=b_2(\Omega)$.  As a consequence
of the Universal Coefficient Theorem, our hypothesis is exactly
equivalent to say that $b_1=0$. By point~(4) above, this is equivalent
to $\chi(\Omega)=1+b_2$ as well. Together with the equality
$\chi(\Omega)= h+1-\sum_{i=0}^h g(S_i)$ proved above, this implies
that
\begin{equation}\label{ind:eq}
h-\sum_{i=0}^h g(S_i)=b_2 \geq 0\ .
\end{equation}

The proof proceeds now by induction on $h\geq 0$. If $h=0$, then we
have $ -g(S_0) \geq 0$, so $g(S_0)=0$ and $S_0$ is a smooth 2--sphere
embedded in $S^3$. Hence, in this case, our theorem reduces to the
celebrated Alexander Theorem (1924) \cite{ALEX} (see also
\cite{HATCH2} for a very accessible proof in the case of smooth
spheres, rather than polyhedral ones as in the original paper by
Alexander). If $h \geq 1$, then equation~(\ref{ind:eq}) implies that
$g(S_{j_0})=0$ for at least one $j_0\in \{0,\ldots,h\}$. Suppose
$j_0\geq 1$. Let us denote by $\Omega^0$ the domain
$\Omega^0=\Omega\cup \overline{\Omega (S_{j_0})}$ obtained by
capping--off the boundary sphere $S_{j_0}$ of $\Omega$ with the
3--disk $\overline {\Omega(S_{j_0})}$. An elementary application of
the Mayer--Vietoris Theorem (see e.g.~\cite{HATCH})
shows that $\Omega^0$ is a domain with
$(h-1)$ boundary components such that $H^1(\Omega^0;\R)=0$, and this
allows us to conclude by induction. If $j_0=0$, then the same proof
applies, after defining $\Omega^0$ as the domain obtained by filling
$\Omega$ (in $S^3$) with the 3--disk $\overline{\Omega^*(S_0)}$.~\cvd

\begin{remark}\label{nec-hyp}
{\rm Theorem~\ref{simple-top} does not hold in general if we don't
assume $\Omega$ to have locally flat boundary.  In fact, on one hand,
the Jordan--Brower Separation Theorem (which is more sophisticated than
Proposition \ref{divide}, see~\cite{brown2}) establishes that every
{\it topological} 2--sphere $S$ embedded in $S^3$ disconnects $S^3$ in
two domains each of which has trivial singular $1$--homology
module. On the other hand, Alexander again (\cite{ALEX2, ALEX3}, see
also \cite[p. 76 and p. 81]{rolfsen}) produced celebrated examples
of non--locally flat topological 2--spheres whose complement in $S^3$
consists of domains one of which (or even both of which) is not simply
connected.}
\end{remark}

\begin{remark}\label{homologysphere}
{\rm A smooth compact connected $3$--manifold $M$ with non--empty
boundary is a \emph{$\Z$--homology disk} (resp.~\emph{$\R$--homology
disk}) if its homology modules with coefficients in $\Z$ (resp.~in
$\R$) are trivial, except that in dimension 0 (so a $\Z$--homology
disk is necessarily a $\R$--homology disk).  Non--simply connected
$\R$--homology disks are easily constructed by removing a small
genuine $3$--ball from closed $3$--manifolds with finite (but
non--trivial) fundamental group such as the projective space
$\mathbb{P}^3(\R)$ or any lens space $L(p,q)$
(see~\cite[p. 233]{rolfsen}).  In the same spirit, a non--simply
connected $\Z$--homology disk can be obtained by removing a genuine
$3$--ball from a closed non--simply connected $3$--manifold having
trivial 1--dimensional $\Z$--homology.  The first example of such a
manifold is due to Poincar\'e.  Theorem~\ref{simple-top} implies that
non--simply connected $\R$--homology disks cannot be embedded in
$S^3$.}
\end{remark}

\begin{remark}
{\rm Even in the locally flat case, the conclusions of
 Theorem~\ref{simple-top} are no longer true when dealing with
 domains in higher dimensional Euclidean space. For example, the
 projective plane $\mathbb{P}^2(\R)$ can be emdedded in $\R^4$, and a
 tubular neighbourhood  of the image of such an embedding  is a
 4--dimensional $\R$--homology disk with fundamental group isomorphic to
 $\Z/2\Z$. }
\end{remark}

We end this section with an open question (as far as we know):
\begin{ques}
{\rm Let $\Omega$ be a {\it not necessarily bounded} domain with
smooth boundary. Assume that $H^1(\Omega;\R)=0$. Does it hold anyway
that $\Omega$ is simply connected?}
\end{ques}

%%%%%%%%%%%%%%%%%%%%%

\section{Helmholtz domains} \label{Helm}

Let us give a definition that covers many current instances in the
literature about Helmholtz cuts method (see also Remark~\ref{rem:very-simple}).

\begin{defi}\label{hel-dom}
{\rm A domain $\Omega \subset \R^3$ with locally flat boundary is {\it
Helmholtz} if there exists a finite family $\Ff=\{\Sigma_i\}$ (called
{\it cut--system for $\Omega$}) of disjoint properly embedded
(connected) surfaces in $(\overline{\Omega},\partial \Omega)$, 
%with non--empty boundary, 
such that every connected component $\Omega^0$ of $\Omega_C(\Ff)$
(i.e.~the disjoint union of domains obtained by cut/open
simultaneouosly along all the $\Sigma_i$'s) satisfies
$H^1(\Omega^0;\R)=0$.}
\end{defi}

We are going to provide an exhaustive and simple characterization of
Helmholtz domains (and of their cut--systems). 
We say that a cut--system for $\Omega$ is \emph{minimal} if it does not properly
contain any cut--system for $\Omega$.
Of course, every cut--system
contains a minimal cut--system.

\begin{lem}\label{helmconnesso}
Suppose $\Ff$ is a minimal cut--system for $\Omega$. Then 
$\Omega_C (\Ff)$ is connected. In particular, every surface
of $\Ff$ has non--empty boundary.
\end{lem}
\begin{proof}
Let
$\Omega_1,\ldots,\Omega_k$ be the connected components of $\Omega_C
(\Ff)$ and suppose by contradiction $k\geq 2$. 
Then we can find a connected surface
$\Sigma_0\in \Ff$ which lies ``between'' two distinct
$\Omega_i$'s. We will now show that the family
$\Ff'=\Ff\setminus \{\Sigma_0\}$ is a cut--system
for $\Omega$, thus obtaining the desired contradiction. 

Up to reordering the $\Omega_i$'s, we may suppose that (parallel
copies of) $\Sigma_0$ lie in the boundary of both $\Omega_{k-1}$ and
$\Omega_k$, so that $\Omega_C (\Ff')=\Omega'_1 \cup \ldots \cup
\Omega'_{k-1}$, where $\Omega'_i=\Omega_i$ for every $i \in
\{1,\ldots,k-2\}$, $\Sigma_0$ is properly embedded in $\Omega'_{k-1}$
and $\Omega_{k-1}\cup \Omega_{k}$ is obtained by cutting
$\Omega'_{k-1}$ along $\Sigma_0$. 
Since $\Ff$ is a cut--system for $\Omega$, the modules
$H^1(\Omega_{k-1};\R)$ and $H^1(\Omega_k;\R)$ are null. By Theorem
\ref{simple-top}, it follows that $\Omega_{k-1}$ and $\Omega_k$ are simply
connected. But $\Sigma_0$ is connected, so an easy application of
Van--Kampen's Theorem (see e.g.~\cite{HATCH}) ensures that
$\Omega'_{k-1}$ is also simply connected, whence $H^1(\Omega'_{k-1};\R)=0$.
Therefore $\Ff'$ is a cut--system for $\Omega$.

We have thus proved the first statement of the lemma. 
Now the conclusion follows from the fact that 
every smooth surface $S\subset \overline\Omega$ without boundary disconnects
$S^3$ (see Proposition~\ref{divide}), whence \emph{a fortiori} $\Omega$.
\end{proof}

\begin{defi}\label{H-B}
{\rm A 3--dimensional {\it $1$--handle} is a $3$--manifold $M$
homeomorphic to $D^2\times [0,1]$ on which there is fixed a
distinguished subspace $A\subset M$ such that the pair $(M,A)$ is
homeomorphic to the pair $(D^2\times [0,1],D^2\times \{0,1\})$. The
connected components of $A$ are the {\it attaching $2$--disks} of $M$,
while if $B\subset M$ corresponds to $D^2\times \{1/2\}$ under a
homeomorphism $(M,A) \cong (D^2\times [0,1],D^2\times \{0,1\})$, then
$B$ is a \emph{co--core} of $M$.  A {\it handlebody} $\overline{H}$ in
$S^3$ is the closure of a domain $H\subset S^3$ with {\it connected}
locally flat boundary (called an {\it open} handlebody), which
decomposes as the disjoint union of $3$--disks (the \emph{$0$--handles
of $\overline{H}$}) together with a disjoint union of 1--handles embedded
in $S^3$ in such a way that the following conditions hold: the
internal part of every $1$--handle is disjoint from the internal part
of every $0$--handle, every attaching $2$--disk of every $1$--handle
lies on the spherical boundary of some $0$--handle, and there are no
further intersections between $0$-- and $1$-- handles (in the smooth
case some ``rounding the corners'' procedure is understood).}
\end{defi}

\begin{remark}\label{h-b-rem}
{\rm It is readily seen that a subset $\overline{H}$ of $S^3$ is a
handlebody if and only if it is equal to a regular neighbourhood of a
finite connected spatial graph $\Gamma$ (i.e. a 1--dimensional
compact connected polyhedron) in $S^3$.  $\Gamma$ is called a {\it spine} of
$\overline{H}$.}
\end{remark}

Every open handlebody $H$ is Helmholtz: a cut--system $\Mm$ for $H$ is
easily contructed by taking one co--core for every $1$--handle of
$\overline{H}$, since in this case the result $H_C (\Mm)$ of cutting $H$
along $\Mm$ is just the family of the internal parts of the
$0$--handles of $\overline{H}$, that are 3--balls. It is not hard to see
that, for suitable subfamilies of these co--cores, the result of
cut/open consists of just one 3--ball.  We will refer to such a
subfamily of co--cores as a {\it minimal system of meridian $2$--disks
for}~$H$. An easy argument using the Euler--Poincar\'e characteristic
shows that the number $g(\overline{H})$ of $2$--disks in a minimal system
of meridian $2$--disks for $H$ equals the genus $g(\partial H)$ of
$\partial H$, and is, in particular, independent from the
handle--decomposition of $\overline{H}$.  We will call $g(\overline{H})$ the
\emph{genus} of $\overline{H}$.  Via ``handle sliding'', it can be easily
shown that two handlebodies are (abstractly) homeomorphic if and only
if they have the same {genus}. Recall that 3--disks are the
handlebodies of genus $0$.

We are now ready to state the main result of this section.  We denote
by $\Omega$ a domain of $\R^3$ with locally flat boundary and by
$S_0,\ldots,S_h$ the connected components of $\partial \Omega$,
ordered as in Lemma~\ref{smoothdom2}.

\begin{teo}\label{helm-char}
$\Omega$ is a Helmholtz domain if and only if the following two
conditions hold:
\begin{itemize}
 \item[$(1)$] The domains $\Omega(S_0)$ and $\Omega^*(S_j)$,
$j=1,\dots,h$, are open handlebodies in $S^3$.
 \item[$(2)$] Every $\Omega(S_j)$, $j=1,\dots, h$, is contained in a
 $3$--disk of $S^3$, embedded in $\Omega(S_0)$, and these $3$--disks
 are pairwise disjoint.
\end{itemize}
Moreover, if $\Omega$ is Helmholtz, then 
there exists 
a cut--system
$\Ff$ for $\Omega$ such that each element of $\Ff$ is a properly
embedded $2$--disk in $(\overline{\Omega},\partial \Omega)$, and
$\Omega_C(\Ff)$ consists of one {``external''} $3$--ball with some
{``internal''} pairwise disjoint $3$--disks removed. In particular,
$\Omega_C(\Ff)$ is connected $($whence simply connected$)$.
\end{teo}

\Dim We can suppose as usual that $\Omega$ has smooth boundary. Assume
that $\Omega$ verifies $(1)$ and $(2)$. Thanks to these conditions, it
is possible to choose a minimal system $\Mm_0$ of meridian $2$--disks
for $\Omega(S_0)$ and, for every $i \in \{1,\ldots,h\}$, a minimal
system $\Mm_i$ of meridian $2$--disks for $\Omega^\ast (S_i)$ in such
a way that $2$--disks belonging to distinct $\Mm_i$'s,
$i=0,1,\ldots,h$, are pairwise disjoint. It is now readily seen that
$\bigcup_{i=0}^h\Mm_i$ provides the cut--system required in the last
statement of the theorem. In particular, $\Omega$ is Helmholtz.

Let us concentrate on the converse implication.  Denote by $\Ff$ an
 arbitrary cut--system for the Helmholtz domain $\Omega$. Accordingly
 to the definition of the cut/open operation along $\Ff$, we have
 $\Omega_C(\Ff) = \Omega \setminus \bigcup_{\Sigma \in \Ff}U_{\Sigma}$,
 where each $U_{\Sigma}$ is a bicollar of $(\Sigma,\partial
 \Sigma)$ in $(\overline{\Omega},\partial \Omega)$, and these
 bicollars are pairwise disjoint. Hence $\overline{\Omega}$ can be
 reconstructed starting from $\overline{\Omega_C(\Ff)}$ by attaching
 to its boundary the $U_{\Sigma}$'s along the surfaces $\Sigma^+$ and
 $\Sigma^-$ corresponding to $\Sigma \times \{\pm 1\}$ in $\Sigma
 \times [-1,1] \cong U_{\Sigma}$. By Theorem~\ref{simple-top}, every
 component of $\Omega_C(\Ff)$ consists of an {``external''} 3--ball
 with some {``internal''} pairwise disjoint 3--disks removed, so the
 boundary components of $\Omega_C(\Ff)$ are spheres. It follows that
 every surface $\Sigma$ is planar, whence homeomorphic either to the
$2$--sphere or to $D^2_k$ for
 some non--negative integer $k$, where $D^2_k$ is
the closure in $\R^2$ of a
 $2$--disk $D^2$ with $k$ disjoint $2$--disks removed from its
 interior.

We will conclude the proof of the theorem in two steps. We will first assume 
that all the surfaces of a given cut--system $\Ff$
of the Helmholtz domain $\Omega$ are $2$--disks. Next we will show how
every arbitrarily given cut--system $\Ff$ can be eventually replaced
with one consisting of $2$--disks only.

{\it Step 1.} Suppose that $\Ff$ is a cut--system for $\Omega$
consisting of $2$--disks only. By Lemma~\ref{helmconnesso},
up to replacing $\Ff$
with a minimal cut--system contained in $\Ff$, we may suppose
that $\Omega_C (\Ff)$ is connected, so that it consists of
 just one ``external'' $3$--ball $B_0$ with some ``internal'' pairwise disjoint
$3$--disks removed. 
%Let
%$\mathcal{B}_1,\ldots,\mathcal{B}_{\ell}$ be the above mentioned
%external $3$--balls. It might happen that $\overline{\mathcal{B}_i}
%\cap \overline{\mathcal{B}_j}\neq \emptyset$, i.e that some of these
%balls are nested.  To overcome this fact, possibly up to considering a
%subfamily of $\Ff$, we can suppose that $\Ff$ is {\it minimal}, that
%is each proper subfamily of $\Ff$ is not a cut--system
%for~$\Omega$. Since $\Omega$ is connected and $\Ff$ is minimal, we
%infer that now $\overline{\mathcal{B}_i} \cap
%\overline{\mathcal{B}_j}=\emptyset \ $ for every $i,j \in
%\{1,\ldots,\ell\}$ with $i \neq j$. 
Observe that we can reconstruct
$\overline{\Omega}$ starting from $\overline{\Omega_C (\Ff)}$ simply
by attaching to $\overline{\Omega_C (\Ff)}$ one $1$--handle for each
$2$--disk in $\Ff$: the attached $1$--handle just coincides with the
removed tubular neighbourhood $D^2\times [0,1]$ of such a $2$--disk in
$\overline{\Omega}$, in such a way that the attaching $2$--disks are
identified with $D^2\times \{0,1\}$.  Let us consider first the
$1$--handles attached to $\overline B_0$. By the very definitions,
%As $\Omega$ is connected, 
the internal part $\Omega (S_0)$ of the union of $\overline B_0$
with such $1$--handles is
%external $3$--balls of $\Omega_C (\Ff)$ is itself connected, and is
%therefore 
an open handlebody. Let $T_1,\ldots,T_h$ be the internal
boundary spheres of $\Omega_C (\Ff)$ and, for each $j \in
\{1,\ldots,h\}$, let $B_j$ be the internal part of the $3$--disk $D_j$
bounded by $T_j$. Now $\overline{\Omega}$ is obtained by attaching to
each $T_j$ some 1--handles contained in the corresponding $D_j$. This
description provides a realization of each $\Omega^*(S_j)$,
$j=1,\dots, h$, as an open handlebody. Note that every $\Omega(S_j)$,
$j=1,\dots, h$, is contained in the corresponding $B_j$. Moreover,
$\Ff$ coincides with the family obtained by taking one co--core
$2$--disk for each added $1$--handle. This completes the proof in the
special case.
\smallskip

{\it Step 2.} Denote by $\Ff$ an arbitrary cut--system for the
Helmholtz domain $\Omega$. Let us show that it is possible to replace
$\Ff$ with a cut--system containing only $2$--disks. 

Up to replacing $\Ff$ with a minimal cut--system, we may assume that
every element of $\Ff$ is homeomorphic to $D^2_k$ for some non--negative
$k$, and that
$\Omega_C (\Ff)$ is connected, so that it consists of
 just one external $3$--ball $B_0$ with some internal pairwise disjoint
$3$--disks $D_1,\ldots, D_l$ removed. We denote by $T_0$ the $2$--sphere
bounding
$B_0$ and by $T_i$ the $2$--sphere bounding $D_i$, $i=1,\ldots,l$, and we 
observe that, under the above assumptions, for every surface
$\Sigma\in\Ff$, there exists $i\in\{0,\ldots,l\}$ 
such that both $\Sigma^+$ and $\Sigma^-$ are contained
in $T_i$. 

\begin{figure}[htbp]
\begin{center}
 \includegraphics[width=9cm]{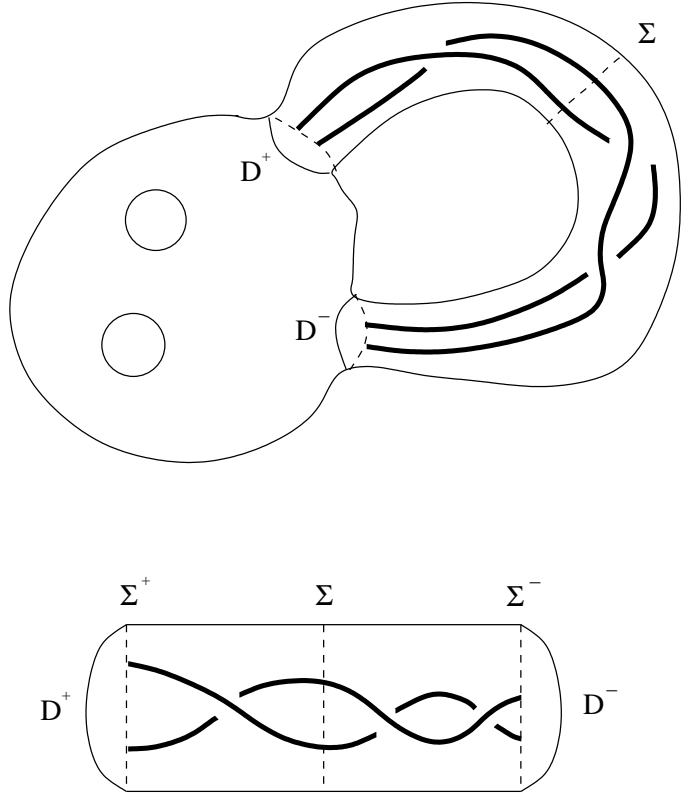}
\caption{\label{trade} Diskal vs planar co--cores: the dashed lines
represent $\Sigma$, $\Sigma^+$ and $\Sigma^-$, while the thickened
strings represent the ``holes'' of $D_k^2\times [-1,1]$ (here $k=2$).}
\end{center}
\end{figure}

We will now show that, if $\Sigma\in\Ff$
is homeomorphic to $D^2_k$ for some $k\geq 1$, then we can obtain a new cut--system
$\Ff'$ from $\Ff$ by replacing $\Sigma$ with two properly embedded $2$--disks.
Such a cut--system will contain a minimal cut--system $\Ff''$ 
with a smaller number (with respect to $\Ff$) of non--diskal surfaces. 
Together with an obvious inductive argument, this will easily imply
that, if $\Omega$ is Helmholtz, then it admits a cut--system consisting of
$2$--disks only, whence the conlusion.
So let $T_i$ be the component of $\partial \Omega_C (\Ff)$ 
containing $\Sigma^+$ and $\Sigma^-$, choose a boundary component $\gamma$
of $D^2_k$ and denote by $\gamma^+$, $\gamma^-$ the curves on $T_i$
corresponding to $\gamma\times\{-1\}$, $\gamma\times\{1\}$ under the identification
of $\Sigma\times\{\pm 1\}$ with $\Sigma^+ \subset T_i$ and $\Sigma^-\subset T_i$.
Now if $D_{{\gamma}^+}$ %(resp.~$D_{{\gamma}^-}$) 
is the $2$--disk on $T_i$
bounded by $\gamma^+$ %(resp.~by $\gamma^-$) 
and containing $\Sigma^+$,% (resp.~$\Sigma^-$), 
we slightly push the internal part of $D_{\gamma^+}$ %(resp.~$D_{\gamma^-}$) 
into $\Omega_C (\Ff)$ thus obtaining
a $2$--disk $D^+$ %(resp.~$D^-$) 
properly embedded in $\Omega_C (\Ff)$
such that $\partial D^+=\gamma^+$ (see Figure~\ref{trade}).
%(resp.~$\partial D^-=\gamma^-$). 
The same procedure applies to $\gamma^-$ providing a $2$--disk $D^-$ properly embedded in $\Omega_C (\Ff)$, and 
of course we may also assume that $D^+$ and $D^-$ are disjoint. Also observe that
by construction both $D^+$ and $D^-$ are disjoint
from every surface in $\Ff$. 

We now set $\Ff'=(\Ff\setminus \{\Sigma\})\cup \{D^+,D^-\}$.
It is easy to see that $\Omega_C (\Ff')$ is given by the disjoint
union of a domain homeomorphic to $\Omega_C (\Ff)$ and a domain
$\Omega'$ 
homeomorphic to the internal part of
$$
\left( D^2 \times [-1,-1+\varepsilon]\right)
\cup \left( D^2_k \times [\varepsilon,1-\varepsilon]\right)
\cup \left( D^2 \times [1-\varepsilon,1]\right) \ .
$$
Now $\Omega'$ is homeomorphic to a $3$--ball with $k$ pairwise
disjoint $3$--disks removed, and is therefore simple
(in the sense of Theorem~\ref{simple-top}). Together
with the fact that $\Omega_C (\Ff)$ is simple, this implies
that $\Ff'$ is a cut--system for $\Omega$.
\cvd

\begin{remark} \label{rem:very-simple}
{\rm 
Bearing in mind the proof of Theorem \ref{helm-char}, we can
now list some equivalent reformulations of the Helmholtz condition
for spatial domains.

$(1)$ 
A domain $\Omega$ of $\R^3$ with locally flat boundary is
Helmholtz if and only if there exists a finite family $\Ss=\{S_i\}$ of
simple domains of $\R^3$ $($in the sense of Theorem
\ref{simple-top}$)$, whose closures are pairwise disjoint, such that
$\overline{\Omega}$ can be constructed starting from the union of the
closures of the $S_i$'s, by attaching some pairwise disjoint
$1$--handles to the boundary spheres of such a union. In addition
$($and equivalently$)$, one may suppose that $\Ss$ consists of a
single simple domain.

$(2)$ A domain $\Omega$ of
$\R^3$ with locally flat boundary is Helmholtz if there exists a
finite family $D_i$, $i=1,\ldots,\ell$, of properly embedded
$2$--disks in $(\overline{\Omega},\partial \Omega)$ such that $\Omega
\setminus \bigcup_{i=1}^{\ell}D_i$ is simply connected. 
In particular, as already mentioned in Lemma~\ref{helmconnesso},
we would get an equivalent definition of Helmholtz domains if we
admitted only cutting surfaces with \emph{non--empty} boundary.

$(3)$ Suppose that $\Omega$ is Helmholtz. 
Then $\Omega$ is \emph{weakly--Helmholtz},
and every cut--system for $\Omega$ is a \emph{weak cut--system} for $\Omega$
(see Section~\ref{weakly-Helmholtz} for the definitions of weakly--Helmholtz domain
and weak cut--system). In particular, Proposition~\ref{charcut}
implies that every cut--system for $\Omega$ contains at least
$b_1 (\Omega)$ surfaces. On the other hand, if 
$\Ff=\{D_1,\ldots, D_\ell\}$ is a cut--system for $\Omega$ consisting
of properly embedded
$2$--disks in $(\overline{\Omega},\partial \Omega)$ such that $\Omega
\setminus \bigcup_{i=1}^{\ell}D_i$ is simply connected, then 
an easy application of the
Mayer--Vietoris Theorem implies that $\ell$ is equal to $b_1(\Omega)$.
Therefore $b_1 (\Omega)$ provides the optimal lower bound on the number
of surfaces contained in the cut--systems for $\Omega$.
}
\end{remark}

\begin{figure}[htbp]
\begin{center}
 \includegraphics[width=9.5cm]{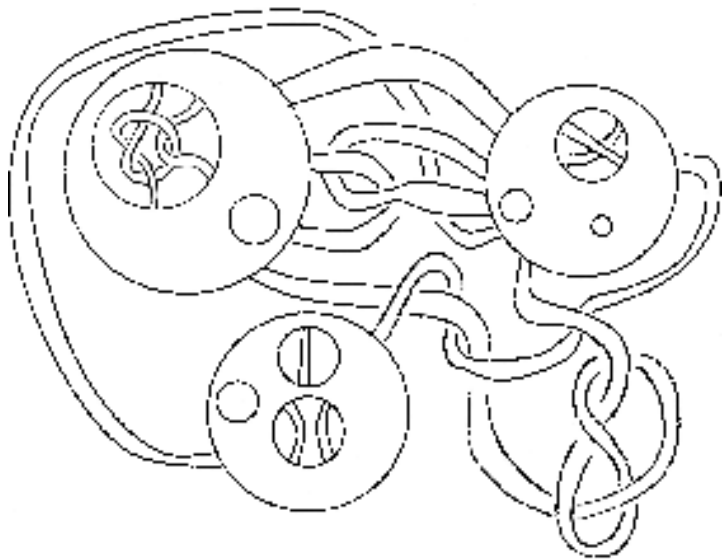}
\caption{\label{helm} A Helmholtz domain.}
\end{center}
\end{figure}

In Figure \ref{helm}, it is drawn a ``typical example'' of
Helmholtz domain: each big circle containing smaller circles
represents an {``external''} $3$--ball with some {``internal''}
pairwise disjoint $3$--disks removed, and the remaining bands
represent the attached $1$--handles.

In some sense, Theorem~\ref{helm-char} should be considered a {\it
negative} result, as it shows that the topology of Helmholtz domains
is forced to be very simple. The following corollary provides an
evidence for this claim. Its proof follows immediately from 
Theorem~\ref{helm-char} and the discussion in Subsection
\ref{link}. For simplicity, we say that a link $L$ of $S^3$ is {\it
Helmholtz} if its complement--domain $\compl(L)$ is.

\begin{cor}\label{helm-link}
Given a link $L$ in $S^3$, the following assertions are equivalent:
\begin{itemize}
 \item[$(1)$] $L$ is Helmholtz.
 \item[$(2)$] $L$ is trivial.
 \item[$(3)$] $\B(L)$ is Helmholtz.
 \end{itemize}
\end{cor}
%\cvd

\begin{figure}[htbp]
\begin{center}
 \includegraphics[width=4.68cm]{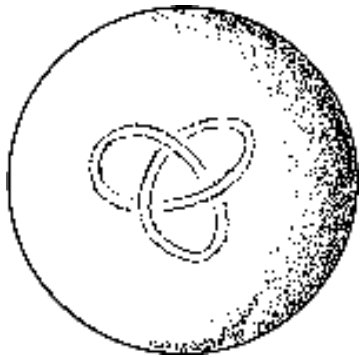}%{box-trefoil-new.eps}
\caption{\label{box-trefoil} A box--domain of a trefoil knot is
not Helmholtz.}
\end{center}
\end{figure}

\noindent The trefoil knot $L$ is not trivial, so the associated
box--domain $\B(L)$, drawn in Figure \ref{box-trefoil}, is a simple
example of a domain of $\R^3$ with smooth boundary, which is not
Helmholtz.

%%%%%%%%%%%%%%%%%%%%%

\subsection{Unknotting reimbedding}\label{reimbedding}
The handlebodies occurring in Theorem \ref{helm-char} are in general
knotted. Let us make precise this notion. A handlebody $\overline{H}$ is
{\it unknotted} if, up to ambient isotopy, it admits a {\it planar}
spine (in the sense of Remark \ref{h-b-rem}) contained in $\R^2
\subset \R^3 \subset S^3$.  Thanks to a celebrated theorem of
Waldhausen \cite{WALD, scharlemann2}, this is equivalent to the fact
that also the complementary domain in $S^3$ is a handlebody: in fact,
a decomposition of $S^3$ into complementary handlebodies is a
so--called {\it Heegaard splitting} of $S^3$, and the Heegaard
splitting of the sphere has been proved to be unique up to isotopy. By
extending the notions of Subsection~\ref{link}, we define a {\it link
of handlebodies} in $S^3$ to be the union of a finite family of
disjoint handlebodies. Such a link is {\it trivial} if all
handlebodies of the family are unknotted and geometrically unlinked,
that is contained in pairwise disjoint 3--disks of $S^3$.

Every (possibly knotted) handlebody can be reimbedded in $S^3$ onto an
unknotted one. We can apply separately this fact to the handlebodies
$\Omega(S_0)$ and $\Omega^*(S_j)$, $j=1,\ldots,h$, of Theorem
\ref{helm-char} and get the following:

\begin{cor} 
A domain $\Omega$ of $\R^3$ with locally flat boundary 
is Helmholtz if and only if it can be reimbedded in $S^3$ onto
a domain $\Omega'$, which is the complement of a trivial link of
handlebodies.
\end{cor}

As an exercise one can see how to realize such an unknotted reimbedding
of the domain of Figure \ref{helm}, just by changing some (over/under)
crossings of the bands representing the $1$--handles. 

By comparing the previous corollary with the following general (and
non--trivial) reimbedding theorem due to Fox \cite{Fox}, we have a
further evidence of the topological simplicity of Helmholtz domains.

\begin{teo}[Fox reimbedding Theorem]\label{Fox}
Every domain $\Omega$ of $S^3$ with locally flat boundary can be
 reimbedded in $S^3$ onto a domain $\Omega'$, which is the complement
 of a link of handlebodies.
 \end{teo}

%%%%%%%%%%%%%%%%%%%%%

\section{Weakly--Helmholtz domain}\label{weakly-Helmholtz}
In this section, we propose and discuss a strictly weaker notion of
``domains that simplify after suitable cuts''. We believe that
the notion we are introducing captures the substance of the philosophy
of Helmholtz cuts, with the advantage of covering a much wider range
of topological models.

In order to save words, from now on, if $M$ is a compact oriented
3--manifold with locally flat boundary, we call \emph{system of
surfaces} in $M$ any finite family $\Ff=\{\Sigma_i\}$ of disjoint
oriented connected surfaces
properly embedded in $M$. We stress that every element of a system of
surfaces is connected and oriented, and that the elements of such a
system are pairwise disjoint. We begin with a definition in the spirit
of condition (b) of Section~\ref{simple} (see also Remark
\ref{rem:w-helm}).

\begin{defi}\label{cut-simple-dom}
{\rm A domain $\Omega$ with locally flat boundary is {\it
weakly--Helmholtz} if it admits a system of surfaces $\Ff$
(called a {\it weak cut--system for $\Omega$}) such that, for every
connected component $\Omega^0$ of $\Omega_C(\Ff)$, the following
condition holds: the restriction to $\Omega^0$ of every curl--free
smooth vector field defined \emph{on the whole of $\Omega$} is the
gradient of a smooth function on~$\Omega^0$.}
\end{defi}

It readily follows from the preceding definition and from Theorem
\ref{helm-char} that every Helmholtz domain is weakly--Helmholtz.

\smallskip

Just as we did in Section~\ref{simple}, let us give some topological
reformulations of the above defi\-nition.  As usual, it is not
restrictive to work in the framework of domains with smooth boundary.
So let $\Omega$ be a domain with smooth boundary, let $\Ff$ be a
system of surfaces in $\overline\Omega$ and let $\Omega_1,\ldots,
\Omega_k$ be the connected components of $\Omega_C(\Ff)$. For
$j \in \{1,\ldots,k\}$, let also $i_j: \Omega_j\lra \Omega$ be the
inclusion. Then $\Ff$ is a weak cut--system for $\Omega$ if and only if one
of the following equivalent conditions hold:

\smallskip

($\beta_1$) {\it For every $j \in \{1,\ldots,k\}$, the image 
of $i_j^*: H^1_{DR}(\Omega)\lra H^1_{DR}(\Omega_j)$ vanishes.}

\smallskip

($\beta_2$) {\it  For every $j \in \{1,\ldots,k\}$, the image of  
$i_j^*: H^1(\Omega;\R)\lra H^1(\Omega_j;\R)$ vanishes.}

\smallskip 

($\beta_3$) {\it For every $j \in \{1,\ldots,k\}$, the image of 
$(i_j)_*: H_1(\Omega_j;\R)\lra H_1(\Omega;\R)$ vanishes.}

\smallskip

The fact that $\Ff$ is a weak cut--system for $\Omega$ if and only if
$(\beta_1)$ holds is a consequence of the canonical isomorphism
between vector fields and 1--forms, the equivalence between
$(\beta_1)$ and $(\beta_2)$ follows from the naturality of de Rham's
isomorphism, and the equivalence between $(\beta_2)$ and $(\beta_3)$
depends on the duality between cohomology and homology.

\subsection{More results about the algebraic topology of domains}
\label{algtop2:subsec}
Before studying weakly--Helmholtz domains, it is convenient to develop
a bit more of information about the algebraic topology of an arbitrary
domain.  While Theorem~\ref{helm-char} provides an exhaustive
description of Helmholtz domains, the classification of
weakly--Helmholtz domains appears to be a quite difficult issue.  In
order to obtain some partial results in this direction, we will use
less elementary (but still ``classical'') tools such as
\emph{relative} homology and \emph{Lefschetz Duality Theorem}.  In
what follows, we will assume that the reader has some familiarity with
such notions and results, which are exhaustively described for
instance in~\cite{HATCH}.  However, in order to preserve as much as
possible the geometric (rather than algebraic) flavour of our
arguments, we will often describe algebraic notions in terms of
geometric ones via an extensive use of transversality. More precisely,
we will often exploit the fact that, if $M$ is a smooth oriented
$n$--dimensional manifold with (possibly empty) boundary $\partial M$,
where $n=2,3$, then every $k$--dimensional (rela\-tive) homology class in
$(M,\partial M)$ with integer coefficients can be geometrically
represented by a smooth oriented closed $k$--manifold (properly)
embedded in $M$. Moreover, the {\it algebraic intersection} 
between a $k$--dimensional and a $(n-k)$--dimensional class (which plays
a fundamental r\^ole in several duality theorems) can be realized
geometrically by taking {transverse} geometric representatives of the
classes involved and counting the intersection points with suitable
signs depending on orientations.
\smallskip 

We now fix a domain $\Omega\subset \R^3$ with smooth boundary. Define
$B:=S^3\setminus \overline \Omega$ and observe that $\partial
\Omega=\overline \Omega \cap \overline B$ is the common smooth
boundary of $\Omega$ and $B$. Since $\partial \Omega$ admits a
bicollar, we may apply the Mayer--Vietoris machinery to the splitting
$S^3=\overline {\Omega} \cup \overline{B}$, obtaining the short exact
sequences
\begin{equation}\label{seconda}
H_2(\partial \Omega;\Z)\lra H_2(\overline \Omega;\Z)\oplus 
H_2(\overline B;\Z)\lra H_2 (S^3;\Z)=0
\end{equation}
and
\begin{equation}\label{prima}
0=H_2 (S^3;\Z)\lra H_1(\partial \Omega;\Z)\lra 
H_1(\overline \Omega;\Z)\oplus H_1(\overline B;\Z)\lra 
H_1 (S^3;\Z)=0 \ .
\end{equation} 

As an immediate consequence, we get the following lemma.

\begin{lem}\label{surj}
The maps $i_\ast\colon H_1 (\partial \Omega;\Z) \lra H_1 (\overline
\Omega;\Z)$ and $i_\ast\colon H_2 (\partial \Omega;\Z)\lra H_2
(\overline \Omega;\Z)$, induced by the inclusion $i \colon \partial
\Omega \hookrightarrow \overline{\Omega}$, are surjective.
\end{lem}

\begin{remark}
{\rm We sketch here a further geometric and more intuitive proof of
the last lemma.
Every class in $H_1(\overline \Omega;\Z)$
can be represented by a knot $C$ embedded in $\Omega$.  Let $S\subset
S^3$ be a Seifert surface for $C$, which we can assume to be
transverse to $\partial \Omega$. Then $S \cap \overline \Omega$
realizes a cobordism between $C$ and a smooth curve contained in
$\partial \Omega$, thus proving that $C$ is homologous to a $1$--cycle
in $\partial \Omega$.

Every class in $H_2(\overline \Omega;\Z)$ can be
represented by the disjoint union of a finite number of compact smooth
orientable surfaces embedded in $\Omega$. 
Every such surface necessarily separates $S^3$ (see
Proposition~\ref{divide}), whence $\Omega$, and is therefore
homologically equivalent to a linear combination of boundary
components.}
\end{remark}

Recall that, if we denote by $T_n (\overline{\Omega})$ the submodule
  of finite--order elements of $H_n (\overline{\Omega},\Z) \cong
  H_n(\Omega;\Z)$, then $T_n (\overline{\Omega})$ is finite for every
  $n \in \mathbb{N}$ and trivial for every $n>2$.

\begin{lem}[see also \cite{GrKo}]\label{free}  
It holds: $T_n (\overline{\Omega})=0$ for every $n\in\mathbb{N}$.
\end{lem}
\Dim Of course, it is sufficient to consider the cases $n=0,1,2$.
Since the $0$--dimensional homology module of any topological space is
free, we have $T_0 (\overline{\Omega})=0$. Moreover, the short exact
sequence~(\ref{prima}) implies that $H_1(\overline{\Omega};\Z)$ is
isomorphic to a submodule of the free $\Z$--module $ H_1(\partial
\Omega;\Z)$, and is therefore free. Finally, by the Lefschetz Duality
Theorem, we have
$$ H^3(\overline \Omega;\Z) \cong H_{0}(\overline \Omega,\partial
\Omega;\Z)=0 \, ,
$$
while the Universal Coefficient Theorem for cohomology gives
$$ H^3(\overline \Omega;\Z) \cong \left( H_3(\overline
\Omega;\Z)/T_3(\overline \Omega)\right) \oplus T_{2}(\overline \Omega)
\, ,
$$
so
$T_2 (\overline{\Omega})=0$.
\cvd

\smallskip

Lemma~\ref{free} implies that the natural morphism
$H_1(\overline{\Omega};\Z)\lra H_1 (\overline{\Omega};\Z) \otimes\R
\cong H_1 (\overline{\Omega};\R)$ is injective. Therefore, keeping
notations from the beginning of Section~\ref{weakly-Helmholtz}, we
obtain that $(\beta_3)$ is equivalent to condition

\smallskip

$(\beta_4)$ {\it For every $j \in \{1,\ldots,k\}$, the image of
$(i_j)_*: H_1(\Omega_j;\Z)\lra H_1(\Omega;\Z)$ vanishes.}

\smallskip

Lemma~\ref{free} allows to describe the Lefschetz Duality Theorem
completely in terms of intersection of cycles. In fact, since $T_0
(\overline \Omega)=0$, the Universal Coefficient Theorem for
cohomology provides a canonical identification $H^1
(\overline{\Omega};\Z) \cong {\rm Hom}(H_1 (\overline{\Omega};\Z),\Z)$
and it turns out that, under the Lefschetz duality isomorphism
$$
H_2(\overline \Omega,\partial \Omega;\Z)\cong 
H^1 (\overline{\Omega};\Z)\cong {\rm Hom}(H_1 (\overline{\Omega};\Z),\Z) \, ,
$$
a class $[\alpha]\in H_2(\overline \Omega,\partial \Omega;\Z)$ is 
identified with the homomorphism which sends every $[\gamma]\in H_1 
(\overline{\Omega};\Z)$ to the algebraic 
intersection between $[\alpha]$ and $[\gamma]$. 
Moreover, since $T_1 (\overline{\Omega})=0$,
a $1$--cycle $\gamma$ in $\overline{\Omega}$ is homologically trivial 
if and only if
its algebraic intersection with every $2$--cycle in 
$H_2 (\overline\Omega,\partial \Omega;\Z)$ is null.

The following lemma will prove useful later.

\begin{lem}\label{geoint}
Let $\Ff=\{\Sigma_1,\ldots, \Sigma_r\}$ be a system of surfaces in
$\overline\Omega$ and let $\gamma$ be a $1$--cycle $($with integer
coefficients$)$ in $\overline{\Omega}$ whose algebraic intersection
with every $\Sigma_i$ is null. Then $\gamma$ is homologous to a
$1$--cycle $\gamma'$ supported in $\overline\Omega\setminus
\bigcup_{i=1}^r \Sigma_i$.
\end{lem}
\Dim Up to homotopy, we may assume that $\gamma$ is the disjoint union
of a finite number of embedded disjoint loops which transversely
intersect $\Sigma_1\cup\ldots\cup \Sigma_r$ in $k$ points $p_1,\ldots,
p_k \in \Omega$. By an obvious induction argument, it is sufficient to prove
that, if $k>0$, then $\gamma$ is homologous to a $1$--cycle $\gamma'$
intersecting $\Sigma_1\cup\ldots\cup \Sigma_r$ in $(k-2)$ points.

Up to reordering the $\Sigma_i$'s, we may assume that $\gamma\cap
\Sigma_1\neq \emptyset$.  Moreover, since the algebraic intersection
between $\gamma$ and $\Sigma_1$ is null, up to reordering the $p_i$'s,
we may suppose that $\gamma\cap \Sigma_1=\{p_j,\, 1\leq j\leq h\}$ for
some $2\leq h\leq k$, and that $\gamma$ intersects $\Sigma_1$ in $p_1$
and $p_2$ with \emph{opposite} orientations.

Let us choose $\epsilon$ small enough in such a way that $\gamma$
intersects the tubular neighbourhood $N_\epsilon (\Sigma_1)$ of
$\Sigma_1$ (in $\overline{\Omega}$) in $h$ small segments
$\gamma_1,\ldots,\gamma_h$ with $p_i \in \gamma_i$ for every
$i$. Since $\Sigma_1$ is connected, if $\alpha$ is a path on
$\Sigma_1$ connecting $p_1$ and $p_2$, then we can define $\gamma'$ by
removing $\gamma_1$, $\gamma_2$ from $\gamma$ and inserting the paths
obtained by pushing $\alpha$ on the boundary components of $N_\epsilon
(\Sigma_1)$ in $\Omega$. Using the fact that $\gamma$
intersects $\Sigma_1$ in $p_1$ and $p_2$ with opposite orientations,
it follows immediately that $\gamma'$ is the disjoint union of a
finite number of embedded loops which can be oriented in such a way
that $[\gamma']=[\gamma]$ in $H_1 (\overline \Omega;\Z)$. This
concludes the proof. \cvd

\smallskip

\noindent
{\bf Assumption:} 
Unless otherwise specified, from now on we only consider homology and
cohomology with \emph{integer} coefficients.

\smallskip

Let us now consider the following portion of the homology exact sequence 
of the pair $(\overline
\Omega,\partial \Omega)$:
\begin{equation}\label{exactseq}
\xymatrix{
H_2(\partial \Omega) \ar[r] &
H_2(\overline \Omega) \ar[r]^{\pi_\ast} &
H_2(\overline \Omega,\partial \Omega) \ar[r]^{\partial} & 
H_1(\partial \Omega) \ar[r]^{i_\ast} & H_1(\overline \Omega)
\ar[r] & H_1(\overline \Omega,\partial \Omega)} \, .
\end{equation}

Let $S_0,\dots,S_h$ be
the boundary components of $\partial \Omega$.

\begin{lem}\label{rA} 
We have the short exact sequence of free modules:
$$
\xymatrix{
0\ar[r] & H_2(\overline \Omega,\partial \Omega)\ar[r]^{\partial} & 
H_1(\partial \Omega)\ar[r]^{\!\!\! i_\ast} 
& H_1(\overline \Omega)\ar[r] & 0 \, .
}
$$
Moreover, ${\rm rank}\, H_2(\overline \Omega,\partial \Omega)=
%{\rm rank}\, H_1(\overline \Omega)=
{\rm rank}\, {\rm Ker}(i_*) =
%{\rm rank}\, {\rm Im}\, \partial = 
b_1 (\overline \Omega)= \sum_{j=0}^{h} g(S_j)$.
\end{lem}
\Dim 
By Lemma~\ref{surj}, the map $\pi_\ast$ in 
sequence~(\ref{exactseq}) is trivial, so 
$\partial$ is injective. Surjectivity of $i_\ast$ and the fact
that $i_\ast\partial =0$ follow respectively by Lemma~\ref{surj} and by
the exactness of sequence~(\ref{exactseq}).
Moreover, we already know that $H_1 (\partial \Omega)$ 
and $H_1 (\overline \Omega)$ are free,
so the sequence splits and $H_2 (\overline{\Omega},\partial \Omega)$ 
is also free.

As a consequence of the exactness of the sequence in the statement, we have
$$
{\rm rank}\, H_2(\overline \Omega,\partial \Omega)=
{\rm rank}\, {\rm Ker}(i_*),\quad
{\rm rank}\, H_1 (\partial \Omega)={\rm rank}\, 
H_2 (\overline \Omega, \partial \Omega) + 
{\rm rank}\, H_1 (\overline \Omega) \, .
$$ Moreover, the Lefschetz Duality Theorem and the Universal
Coefficient Theorem give the isomorphisms $H_2
(\overline{\Omega},\partial \Omega)\cong H^1 (\overline{\Omega})\cong
H_1 (\overline{\Omega})$, so ${\rm rank}\, H_2
(\overline{\Omega},\partial \Omega)=b_1 (\overline{\Omega})$ and hence
${\rm rank}\, H_1 (\partial \Omega)=2 \, {\rm rank}\, H_1 (\overline
\Omega)$, i.e.~$b_1(\partial \Omega)=2b_1 (\overline\Omega)$. But
homology is additive with respect to disjoint union of topological
spaces, so $b_1 (\partial \Omega)= 2\sum_{j=0}^{h} g(S_j)$, whence the
conclusion. \cvd

\begin{remark}\label{uniqueslope}
{\rm Let $K\subset S^3$ be a knot with complement--domain
$\Omega=\compl (K)$.  Lemma~\ref{rA} implies that the kernel of the
map $i_\ast\colon H_1 (\partial\Omega)\to H_1 (\overline\Omega)$ is
freely generated by the class $[\gamma]$ of a non--trivial loop on
$\partial\Omega$.  Let $S$ be a Seifert surface for $K$ intersecting
$\partial \compl (K)$ in a simple loop $\alpha$ parallel to $K$. Since
$\alpha$ bounds the surface $S\cap \overline\Omega$ properly embedded
in $\Omega$, the class $[\alpha]$ is a multiple of $[\gamma]$, and
using that $\alpha$ is simple and not homologically trivial it is not
difficult to show that in fact $[\alpha]=\pm [\gamma]$.  Finally, two
simple closed loops on a torus define the same homology class if and
only if they are isotopic, so we can conclude that the isotopy class
of the loop obtained as the transverse intersection of
$\partial\Omega$ with a Seifert surface for $K$ does not depend on the
chosen surface, as claimed in Example~\ref{Seifert}.}
\end{remark}

\smallskip

\subsection{Cut number and corank}  
It turns out that the property of being weakly--Helmholtz admits
characterizations in terms of classical properties of manifolds and of
their fundamental group. We begin with the following definitions,
which in the case of closed manifolds date back to~\cite{Stallings}
(see also~\cite{Harvey} and~\cite{Sikora}).

\begin{defi}\label{cutnumber}
{\rm Let $M$ be a (possibly non--orientable) smooth connected compact
$3$--manifold with (possibly empty) boundary. The \emph{cut number}
$c(M)$ of $M$ is the maximal number of disjoint properly embedded
(bicollared connected) surfaces $\Sigma_1,\ldots,\Sigma_k$ in
$(M,\partial M)$ such that $M \setminus \bigcup_{i=1}^k \Sigma_i$ is
connected.}
\end{defi}

\begin{defi}\label{corank}
{\rm For each non--negative integer $r$, we denote by $\fr$ the
$r^{\rm th}$--free power of~$\Z$. Given a group $\Gamma$, the
\emph{corank} of $\Gamma$ is the maximal non--negative integer $r$
such that $\fr$ is isomorphic to a quotient of $\Gamma$.}
\end{defi}

Let $M$ be as in Definition~\ref{cutnumber}. It is not difficult to
show that ${\rm corank}\, (\pi_1 (M))\leq b_1 (M)$ (see
e.g.~Corollary~\ref{ovvio} and Remark~\ref{cutcorank:rem}). It was
first observed by Stallings that $c(M)={\rm corank}\, (\pi_1
(M))$. For the sake of completeness, in Proposition~\ref{cutcorank}
below, we will give a proof of such an equality in the case we are
interested in, i.e.~when $M=\overline\Omega$ for some domain $\Omega$
with smooth boundary. Our proof of Proposition~\ref{cutcorank} follows
closely Stallings' original proof (see also~\cite{Sikora}) and can
therefore be easily adapted to deal with the general case.

Before going on, we recall that the elements $a_1,\ldots,a_r$ of a
$\Z$--module $A$ are said to be \emph{linearly independent} if
whenever $c_1,\ldots,c_r\in \Z$ are such that $\sum_{i=1}^r c_i a_i
=0$, then $c_i=0$ for every $i$ (in particular, a set of linearly
independent elements do not contain torsion elements).  We say that a
finite set $a_1,\ldots,a_r$ is a \emph{basis} of $A$ if, for every $a
\in A$, there exists a unique $r$--uple of coefficients
$(c_1,\ldots,c_r)\in\Z^r$ such that $a=\sum_{i=1}^r c_i a_i$ or,
equivalently, if the $a_i$'s are linearly independent and generate
$A$. Of course, if $A$ admits a basis $a_1,\ldots, a_r$, then $A$ is
free of rank $r$. A submodule $\Lambda$ of $A$ is \emph{full} if it is
not a proper finite--index submodule of any other submodule of
$A$. Recall that, if $\Lambda$ is a submodule of $A$, then $\Lambda$
has finite--index in $A$ if and only if ${\rm rank}\, \Lambda={\rm rank}\, A$.
Therefore, if $\Lambda$ is full and ${\rm rank}\, \Lambda={\rm rank}\,
A$, then $\Lambda=A$.

\smallskip

Let now $\Omega\subset \R^3$ be a domain with smooth boundary.  We
define $d(\overline \Omega)$ as the maximal number of \emph{disjoint}
oriented connected surfaces with non--empty boundary, properly embedded
in $(\overline\Omega,\partial\Omega)$, which define linearly
independent elements in $H_2 (\overline{\Omega},\partial \Omega)$.

We begin with the following result.

\begin{lem}\label{linindip}
Let $\Ff=\{\Sigma_1,\ldots,\Sigma_r\}$ be a system of surfaces in
$\overline{\Omega}$ and let $[\Sigma_i]\in H_2( \overline
\Omega,\partial\Omega)$ be the class represented by $\Sigma_i$,
$i=1,\ldots,r$. Then the following conditions are equivalent:
\begin{itemize}
 \item[$(1)$] The $[\Sigma_i]$'s are linearly independent in $H_2
 (\overline \Omega,\partial\Omega)$.
 \item[$(2)$] The $[\Sigma_i]$'s are linearly independent and generate a
 full submodule of $H_2 (\overline \Omega,\partial\Omega)$.
 \item[$(3)$] The set $\Omega_C (\Ff)$ is connected.
\end{itemize}
\end{lem}
\Dim $(1) \Longrightarrow (3)$ Let $\Omega':=\overline\Omega \setminus
\bigcup_{i=1}^r \Sigma_i$. Since $\overline{\Omega_C(\Ff)}$ is a
strong deformation retract of $\Omega'$, it is sufficient to show that
$\Omega'$ is connected. Suppose by contradiction that $\Omega'$ is
disconnected and let $\Omega^0$ be a connected component of $\Omega'$
with $\partial \overline{\Omega^0} \setminus \partial
\Omega=(\Sigma_{j_1} \cup \ldots \cup \Sigma_{j_l}) \setminus \partial
\Omega$ (where $j_h \neq j_k$ if $h\neq k$). Then
$[\Sigma_{j_1}]+\ldots +[\Sigma_{j_l}]=0$ in $H_2
(\overline\Omega,\partial \Omega)$, a contradiction.

$(3) \Longrightarrow (2)$ Recall that, under the Lefschetz duality
isomorphism
$$
H_2(\overline \Omega,\partial \Omega)\cong 
H^1 (\overline{\Omega})\cong {\rm Hom}(H_1 (\overline{\Omega}),\Z) \, ,
$$ the class $[\Sigma_j]\in H_2(\overline \Omega,\partial \Omega)$ is
identified with the linear map $f_j\colon
H_1(\overline{\Omega})\lra\Z$ which sends every $[\gamma]\in H_1
(\overline{\Omega})$ to the algebraic intersection between $\Sigma_j$
and $\gamma$. Now, since $\Omega_C (\Ff)$ is connected, for every $i
\in \{1,\ldots,r\}$, we can construct a loop $\gamma_i \subset \Omega$
which intersects $\Sigma_i$ transversely in one point and is disjoint
from $\Sigma_j$ for every $j\neq i$. It readily follows that, if
$\sum_{j=1}^r c_j f_j=0$, then, for every $i\in \{1,\ldots,r\}$, we
have that $c_i=(\sum_{j=1}^r c_j f_j)(\gamma_i)=0$, so the
$[\Sigma_i]$'s are linearly independent. Let now $\Lambda$ be the
submodule of ${\rm Hom}(H_1 (\overline{\Omega}),\Z)$ generated by the
$f_j$'s and suppose that $\Lambda'$ is a submodule of ${\rm Hom}(H_1
(\overline{\Omega}),\Z)$ with $\Lambda\subset \Lambda'$. Also suppose
that $\Lambda$ has finite--index in $\Lambda'$, and take an element
$f\in\Lambda'$. Our assumptions imply that there exists
$n\in\Z\setminus \{0\}$ such that $n\cdot f$ lies in $\Lambda$ and is
therefore a linear combination $\sum_{i=1}^r c_i f_i$ of the
$f_i$'s. For every $i \in \{1,\ldots,r\}$, it follows therefore that
$c_i=n f(\gamma_i)$, so $c_i=n c'_i$ for some $c'_i\in\Z$ and
$f=\sum_{i=1}^r c'_i f_i\in\Lambda$. We have thus proved that
$\Lambda$ is full.

$(2) \Longrightarrow (1)$ is obvious. \cvd

\smallskip

The following proposition relates to each other the notions just
introduced.

\begin{prop}\label{cutcorank}
It holds:
$$
d(\overline{\Omega})=c(\overline\Omega)=
{\rm corank}\, (\pi_1 (\overline \Omega)) \, .
$$
\end{prop}
\Dim The equality $d(\overline\Omega)=c(\overline\Omega)$ is an
immediate consequence of Lemma~\ref{linindip}. In order to prove the
proposition, we will now prove the inequalities
$c(\overline\Omega)\leq {\rm corank} (\pi_1(\overline{\Omega})) \leq
d(\overline\Omega)$.

So let $\Ff=\{\Sigma_1,\ldots,\Sigma_r\}$ be a system of surfaces in
$\overline\Omega$ such that $\overline\Omega \setminus \bigcup_{i=1}^r
\Sigma_i$ is connected and let $B_r$ be the wedge of $r$ copies
$S^1_1,\dots, S^1_r$ of the unitary circle, with base point $x_0$.
Also recall that the fundamental group $\pi_1(B_r,x_0)$ is freely
generated by the (classes of the) loops $\gamma_1,\ldots,\gamma_r$,
where $\gamma_j\colon [0,1]\lra S^1_j$ is a generator of
$\pi_1(S^1_j,x_0)$ (in particular, $\gamma(0)=\gamma(1)=x_0$).  By a
classical Pontryagin--Thom construction (see \cite{MIL}), we can construct a
continuous map
$$
f=f_\Ff \colon \overline\Omega \lra B_r
$$ as follows. Consider a system of disjoint closed bicollars $U_j$ of
the $\Sigma_j$'s in $\overline{\Omega}$ 
and fix diffeomorphic identifications $U_j\cong
\Sigma_j\times [0,1]$, $j=1,\ldots,r$. Then, if $(x,t)\in U_j$, we set
$f(x,t)=\gamma_j (t)$, while, for $q \in M\setminus \bigcup_{j=1}^r
U_j$, we set $f(q)=x_0$. Since $\overline\Omega \setminus
\bigcup_{j=1}^r U_j$ is connected, it is easily seen that, if $p$ is
any basepoint in $\overline\Omega \setminus \bigcup_{j=1}^r U_j$, then
the map $f_\ast\colon \pi_1 (\overline\Omega,p)\lra \pi_1 (B_r,x_0)$
is surjective. We have thus shown that $c(\overline\Omega) \leq {\rm
corank}(\pi_1 (\overline\Omega))$.

In order to prove that ${\rm corank}(\pi_1 (\overline{\Omega})) \leq
d(\overline\Omega)$, we can invert the construction just described as
follows. Let $r={\rm corank}(\pi_1 (\overline{\Omega}))$ and take a
surjective homomorphism $\phi:\pi_1(\overline \Omega)\lra \fr$. As
$B_r$ is a $K(\fr,1)$ space with contractible universal covering (see
\cite{HATCH}), there exists a continuous surjective map $f: \overline
\Omega \lra B_r$ such that $\phi = f_*$. Up to homotopy, we can assume
that the restriction of $f$ to $f^{-1}(B_r\setminus \{x_0\})$ is
smooth. By the Morse--Sard Theorem (see \cite{MIL,HIR}), we can select
a regular value $x_j\in S^1_j \setminus \{x_0\}$ and define
$N_j:=f^{-1} (x_j)$ for every $j \in \{1,\ldots,r\}$. Then $N_j$ is a
finite union of disjoint properly emdedded surfaces in $\overline
\Omega$. Moreover, if we fix an orientation on every $S^1_j$, then we
can define an orientation on $N_j$ by the usual ``first the outgoing
normal vector'' rule, where a vector $v$ is outgoing in $q\in N_j$ if
$df (v)$ is positively oriented as a vector of the tangent space to
$S^1_j$ in $f(q)$. Let now $p$ be a basepoint in $f^{-1} (x_0)\subset
\overline\Omega$ and let $\alpha_j$ be a loop in $\Omega$ based at $p$
whose homotopy class $[\alpha_j]\in \pi_1 (\overline\Omega,p)$ is
mapped by $\phi=f_\ast$ onto a generator of $\pi_1 (S^1_j,x_0)< \pi_1
(B_r,x_0)$. Up to homotopy, we may suppose that the intersection
between $\alpha_j$ and $N_j$ is transverse. Moreover, by the very
construction of $\alpha_j$, the algebraic intersection between
$\alpha_j$ and $N_k$ is equal to $1$ if $j=k$ and to $0$ otherwise. In
particular, there exists a connected component $\Sigma_j$ of $N_j$
such that the algebraic intersection of $\alpha_j$ with $\Sigma_k$ is
not null if and only if $k\neq j$. By Lefschetz Duality Theorem, this
readily implies that $\Sigma_1,\ldots,\Sigma_r$ represent linearly
independent elements of $H_2(\overline\Omega,\partial\Omega)$. This
gives in turn the inequality ${\rm corank}\, (\pi_1
(\overline{\Omega}))\leq d(\overline\Omega)$. \cvd

\smallskip

Since $d(\overline{\Omega})\leq {\rm rank}\,
H_2(\overline\Omega,\partial \Omega)={\rm rank}\, H_1
(\overline\Omega)$, Proposition~\ref{cutcorank} immediately implies
the following result.

\begin{cor}\label{ovvio}
It holds: $c(\overline\Omega)={\rm corank}\,
(\pi_1(\overline\Omega))\leq b_1(\overline\Omega)$.
\end{cor}

\begin{remark}\label{cutcorank:rem}
{\rm As mentioned above, the relations
$c(M)={\rm corank}\, (\pi_1 (M))\leq b_1 (M)$ hold in general,
i.e.~even when $M$ is any (possibly non--orientable) manifold.  In
fact, the proof of Proposition~\ref{cutcorank} can be easily adapted
to show that $c(M)={\rm corank}\, (\pi_1 (M))$. Moreover, if ${\rm
corank}\, (\pi_1 (M))=r$, then there exists a surjective homomorphism
from $\pi_1 (M)$ to the Abelian group $\mathbb{Z}^r$. As a consequence
of the classical Hurewicz Theorem (see e.g.~\cite{HATCH}), such a
homomorphism factors through $H_1 (M)$, whose rank is therefore at
least $r$.  This readily implies the inequality ${\rm corank}\, (\pi_1
(M))\leq b_1 (M)$.}
\end{remark}

\subsection{Topological characterizations of weakly--Helmholtz domains}
The following lem\-ma shows that, just as in the case of Helmholtz
domains, every weakly--Helmholtz domain admits a non--disconnecting
cut--system. So let $\Omega\subset \R^3$ be a domain with smooth
boundary.

\begin{lem}\label{nosconnesso}
If $\Omega$ is weakly--Helmholtz, then it admits a weak cut--system
whose surfaces do not disconnect $\Omega$. More precisely, every weak
cut--system $\Ff$ for $\Omega$ contains a weak cut--system $\Ff'$ for
$\Omega$ such that $\Omega_C (\Ff')$ is connected.
\end{lem}
\Dim Let $\Ff$ be a weak cut--system for $\Omega$, let
$\Omega_1,\ldots,\Omega_k$ be the connected components of $\Omega_C
(\Ff)$ and suppose $k\geq 2$. Then we can find a connected surface
$\Sigma_0\in \Ff$ which lies ``between'' two distinct
$\Omega_i$'s. Let us set $\Ff'=\Ff\setminus \{\Sigma_0\}$ and show
that $\Ff'$ is a weak cut--system for $\Omega$. By repeating this
procedure $k-1$ times, we will be left with the desired weak
cut--system that does not disconnect $\Omega$.

Up to reordering the $\Omega_i$'s, we may suppose that (parallel
copies of) $\Sigma_0$ lie in the boundary of both $\Omega_{k-1}$ and
$\Omega_k$, so that $\Omega_C (\Ff')=\Omega'_1 \cup \ldots \cup
\Omega'_{k-1}$, where $\Omega'_i=\Omega_i$ for every $i \in
\{1,\ldots,k-2\}$, $\Sigma_0$ is properly embedded in $\Omega'_{k-1}$
and $\Omega_{k-1}\cup \Omega_{k}$ is obtained by cutting
$\Omega'_{k-1}$ along $\Sigma_0$. We now claim that every $1$--cycle
in $\Omega'_{k-1}$ decomposes, up to boundaries, as the sum of a
$1$--cycle supported on $\Omega_{k-1}$ and a cycle supported in
$\Omega_k$. In fact, since $\Sigma_0$ disconnects $\Omega'_{k-1}$, the
homology class represented by $\Sigma_0$ in $H_2
(\Omega'_{k-1},\partial \Omega'_{k-1})$ is null.  This implies that
the algebraic intersection between $\Sigma_0$ and any $1$--cycle in
$\Omega_{k-1}'$ is null, and the claim now follows from
Lemma~\ref{geoint}.

The claim just proved implies that the image of
$(i'_{k-1})_\ast:H_1(\Omega_{k-1}') \lra H_1(\Omega)$ equals the sum
of the images of $(i_{k-1})_\ast \colon H_1(\Omega_{k-1})\lra H_1
(\Omega)$ and of $(i_{k})_\ast \colon H_1(\Omega_{k})\lra H_1
(\Omega)$, which are both trivial, because of $\Ff$ satisfies
condition ($\beta_4$). Therefore the image of $(i'_j)_\ast$ vanishes
for every $j \in \{1,\ldots,k-1\}$, so $\Ff'$ is a weak cut--system
for $\Omega$. \cvd

\smallskip

\begin{lem}\label{generano}
Let $\Ff=\{\Sigma_1,\ldots, \Sigma_r\}$ be a system of surfaces in
$\overline\Omega$ and let $\Lambda \subset
H_2(\overline\Omega,\partial\Omega)$ be the submodule generated by the
classes $[\Sigma_1],\ldots,[\Sigma_r]$ represented by the
$\Sigma_i$'s. The system $\Ff$ is a weak cut--system if and only if
${\rm rank}\, \Lambda=b_1 (\overline\Omega)$.
\end{lem}
\Dim We claim that $\Ff$ is a weak cut--system for $\Omega$ if and only if the 
following condition holds:
\begin{itemize}
 \item if $[\gamma]\in H_1(\overline\Omega)$ has null algebraic
intersection with every $[\Sigma_i]$, $i=1,\ldots,r$, then
$[\gamma]=0$ in $H_1(\overline\Omega)$.
\end{itemize}
In fact, suppose $\Ff$ is a weak cut-system and let 
$[\gamma]\in H_1(\overline\Omega)$ have
null algebraic
intersection with every $[\Sigma_i]$, $i=1,\ldots,r$.
Then,
by Lemma~\ref{geoint}, we can suppose that $[\gamma]$ is represented
by a $1$--cycle supported in $\Omega_C (\Ff)$. This implies that, if
$\Omega_1,\ldots,\Omega_k$ are the connected components of
$\Omega_C(\Ff)$, then $[\gamma]=\sum_{i=1}^k [\gamma_i]$ in $H_1
(\overline\Omega)$, where the $1$--cycle $\gamma_i$ is supported in
$\Omega_i$ for every $i$. But, by condition~$(\beta_4)$, if $\Ff$ is a
weak cut--system, we have $[\gamma_i]=0$ in $H_1(\overline\Omega)$ for
every $i$, so $[\gamma]$ is homologically trivial in $\Omega$. On the
other hand, if the inclusion $i_j\colon \Omega_j\lra\overline\Omega$
induces a non--trivial homomorphism $(i_j)_\ast\colon
H_1(\Omega_j)\lra H_1(\overline \Omega)$, then every non--null class
$[\gamma]$ in ${\rm Im}\, (i_j)_\ast$ has null algebraic intersection
with every $[\Sigma_i]$, $i=1,\ldots,r$. This concludes the proof of
the claim.

For every $j \in \{1,\ldots,r\}$, let now $f_j\colon
H_1(\overline\Omega)\lra\Z$ be the linear map corresponding to
$[\Sigma_j]$ under the identification
$$
H_2 (\overline\Omega,\partial\Omega)\cong {\rm Hom}\, 
(H_1(\overline\Omega),\Z) \, .
$$ The claim above shows that $\Ff$ is a weak cut--system for $\Omega$
if and only if
$$
\bigcap_{i=1}^r {\rm Ker} (f_i) =\{0\} \, .
$$ It is now a standard fact of Linear Algebra that this last
condition is satisfied if and only if the $f_i$'s generate a
finite--index submodule of ${\rm Hom}\, (H_1(\overline\Omega),\Z)$,
whence the conclusion. \cvd

\begin{cor}\label{limitazione:cor}
Every weak cut--system for $\Omega$ contains at least
$b_1(\overline\Omega)$ surfaces.
\end{cor}

We can now summarize the results obtained so far in the following
Proposition~\ref{charcut} and Theorem~\ref{char1}, which provide a
characterization of weakly--Helmholtz domains and of their weak
cut--systems. We begin with the following definition.

\begin{defi}\label{minimal}
{\rm A weak cut--system $\Ff$ for $\Omega$ is \emph{minimal} if every
proper subset of $\Ff$ is \emph{not} a weak cut--system for $\Omega$.}
\end{defi}

It follows by the definitions that every system of surfaces
containing a weak cut--system is itself a weak cut--system, so a
system of surfaces is a weak cut--system if and only if it contains a
minimal weak cut--system.

\begin{prop}\label{charcut}
Let $\Ff=\{\Sigma_1,\ldots,\Sigma_r\}$ be a system of surfaces in
$\overline{\Omega}$, and let $[\Sigma_i]\in
H_2(\overline\Omega,\partial\Omega)$ be the class represented by
$\Sigma_i$, $i=1,\ldots,r$. Then the following conditions are
equivalent.
\begin{itemize}
 \item[$(1)$] $\Ff$ is a minimal weak cut--system for $\Omega$.
 \item[$(2)$] $r=b_1(\overline\Omega)$ and $\Omega_C (\Ff)$ is connected.
 \item[$(3)$] The $[\Sigma_i]$'s provide a basis of
 $H_2(\overline\Omega,\partial\Omega)$.
 \item[$(4)$] $r=b_1 (\overline\Omega)$ and the $[\Sigma_i]$'s are
 linearly independent elements in
 $H_2(\overline\Omega,\partial\Omega)$.
\end{itemize}
\end{prop}
\Dim Let us denote by $\Lambda$ the submodule of $H_2(\overline
\Omega,\partial \Omega)$ generated by the $[\Sigma_{i}]$'s.

$(1) \Longrightarrow (2)$ By Lemma~\ref{nosconnesso}, the minimality
of $\Ff$ implies that $\Omega_C (\Ff)$ is connected. Moreover, by
Lemmas~\ref{linindip} and~\ref{generano}, $\Lambda$ is freely
generated by the $[\Sigma_{i}]$'s and $r=b_1 (\overline{\Omega})$.

$(2) \Longrightarrow (3)$ By Lemma~\ref{linindip}, since $\Omega_C
(\Ff)$ is connected, $\Lambda$ is full and freely generated by the
$[\Sigma_i]$'s. The assumption $r=b_1(\overline\Omega)={\rm rank}\,
H_2(\overline \Omega,\partial \Omega)$ easily implies that $\Lambda$
has finite--index in $H_2(\overline \Omega,\partial \Omega)$. Being
full, $\Lambda$ is then equal to the whole $H_2(\overline
\Omega,\partial \Omega)$, and the $[\Sigma_i]$'s provide therefore a
basis of $H_2(\overline \Omega,\partial \Omega)$.

$(3) \Longrightarrow (4)$ is obvious.

$(4) \Longrightarrow (1)$ Condition~(4) readily implies that ${\rm
rank}\, \Lambda ={\rm rank}\, H_2(\overline \Omega,\partial \Omega)$,
so $\Lambda$ has finite--index in $H_2(\overline \Omega,\partial
\Omega)$. Thanks to Lemma~\ref{generano}, $\Ff$ is a weak cut--system
for $\Omega$. Moreover, $\Ff$ is minimal by
Corollary~\ref{limitazione:cor}. \cvd

\smallskip

As a consequence of Propositions~\ref{cutcorank} and~\ref{charcut}, we
obtain the following characterization of weakly--Helmholtz domains.

\begin{teo}\label{char1}
Let $\Omega\subset \R^3$ be a domain with locally flat boundary and let
$r:=b_1(\overline{\Omega})$. Then the following conditions are
equivalent:
\begin{itemize}
 \item[$(1)$] $\Omega$ is weakly--Helmholtz.
 \item[$(2)$] There exists a system of surfaces
 $\Ff=\{\Sigma_1,\dots,\Sigma_r\}$ in $\overline{\Omega}$ such that
 $\overline{\Omega} \setminus \bigcup_{i=1}^r\Sigma_i$ is connected.
 \item[$(3)$] There exists a basis of $H_2(\overline \Omega,\partial
 \Omega)$ represented by a system of surfaces in $\overline{\Omega}$.
 \item[$(4)$] $c(\overline \Omega)=d(\overline\Omega)={\rm corank}\, 
(\pi_1 (\overline\Omega))=r$.
 \item[$(5)$] There exists a surjective homomorphism from
 $\pi_1(\Omega)$ onto $\fr$.
\end{itemize}
\end{teo}
%\cvd

\begin{remark} \label{rem:w-helm}
{\rm $(1)$ By the preceding theorem, it is possible to give an
equivalent definition of weakly--Helmholtz domain as follows: ``a
domain $\Omega$ of $\R^3$ is weakly--Helmholtz if there exists a
finite family $\{\Sigma_i\}$ of disjoint properly embedded (connected)
surfaces in $(\overline{\Omega},\partial \Omega)$, with non--empty
boundary, such that $\Omega^*:=\Omega \setminus \bigcup_i \Sigma_i$ is
connected and the restriction to $\Omega^*$ of every curl--free smooth
vector field defined on the whole of $\Omega$ is the gradient of a
smooth function on~$\Omega^*$''.}

{\rm $(2)$ As in the case of Helmholtz domains, one can obtain other
equivalent definitions of weakly--Helmholtz domain starting from
Definition \ref{cut-simple-dom} or from the definition given in the
preceding point $(1)$ by admitting only cutting surfaces with
non--empty boundary.}
\end{remark}

Let $L$ be a link in $S^3$. We say that $L$ is {\it weakly--Helmholtz}
if the complement--domain $\compl(L)$ of $L$ is (see Subsection
\ref{link}). We have the following easy:

\begin{lem}\label{equivalence}
The link $L$ is weakly--Helmholtz if and only if its box--domain $\B (L)$
is. 
\end{lem}
\Dim Recall that $\B (L)$ is obtained by removing a small $3$--disk
$D$ from $\compl (L)$. An easy application of the Mayer--Vietoris
machinery now implies that the modules $H_1 ({\compl (L)})$
and $H_1 ({\B (L)})$ are isomorphic, so $b_1
({\compl (L)})=b_1 ({\B (L)})$. On the other hand,
an easy application of Van Kampen's Theorem (see e.g.~\cite{HATCH})
ensures that the fundamental groups $\pi_1 ({\compl (L)})$
and $\pi_1 ({\B (L)})$ are also isomorphic, so $b_1
({\compl (L)})={\rm corank}\, (\pi_1 ({\compl
(L)}))$ if and only if $b_1 ({\B (L)})={\rm corank}\, (\pi_1
({B (L)}))$. Now the conclusion follows from
Theorem~\ref{char1}.  \cvd

\smallskip 

As a consequence of Corollary~\ref{helm-link}, we know
that a knot in $S^3$ is Helmholtz if and only if it is trivial. On the
contrary, every knot is weakly--Helmholtz as we see in the next
result.

\begin{cor}\label{knots}
The following statements hold.
\begin{itemize}
 \item[$(1)$] Every knot in $S^3$ is weakly--Helmholtz.
 \item[$(2)$] The box--domain of any knot in $S^3$ is weakly--Helmholtz.
 \end{itemize}
\end{cor}
\Dim Let $S$ be a Seifert surface of a knot $K$ in $S^3$. Since $S$
does not disconnect the complement--domain $\compl(K)$ of $K$, the
equivalence $(1) \Longleftrightarrow (4)$ in Theorem~\ref{char1}
immediately implies that $K$ is weakly--Helmholtz. Therefore $(1)$ is
proved, and  $(2)$ now follows from Lemma~\ref{equivalence}. 
%The box--domain $\B(K)$ of $K$ can be
%identified with $\compl(K)$ with a small $3$--disk $D$ removed, so, up
%to a slight deformation of $S$, we may suppose that $S$ does not
%intersect $D$. Evidently, $\B(K) \setminus S$ is connected, and hence
%$(2)$ follows by applying Proposition~\ref{char1} again. 
\cvd

\begin{remark}
{\rm The box--domain of a trefoil knot, drawn in above Figure
\ref{box-trefoil}, is a simple example of weakly--Helmholtz, but not
Helmholtz, domain.}
\end{remark}

%%%%%%%

\subsection
{The intersection form on surfaces} Let $S$ be a connected compact
orientable surface. If $\alpha,\beta$ are $1$--cycles on $S$, up to
homotopy, we can suppose that $\alpha$ and $\beta$ transversely
intersect in a finite number of points, and define the algebraic
intersection between $\alpha$ and $\beta$ as the difference between
the number of points in which they intersect ``positively'' and the
number of points in which they intersect ``negatively'', with respect
to the fixed orientation on $S$. It is not difficult to show that the
algebraic intersection defines a bilinear \emph{skew--symmetric}
product on the space of $1$--cycles, and that the algebraic
intersection between a boundary and any $1$--cycle is null. It follows
that such a bilinear product descends to homology, thus defining a
bilinear skew--symmetric intersection form
$$ \langle\ \cdot\ , \ \cdot \ \rangle\colon H_1 (S)\times H_1(S)\lra
\Z \, .
$$ Being a particular instance of the general Lefschetz Duality
Theorem just recalled, such an intersection form induces an
isomorphism between $H_1(S)$ and ${\rm Hom}\, (H_1(S),\Z)\cong
H^1(S)$. In particular, $H_1(S)$ admits a \emph{symplectic basis},
i.e.~a free basis $\alpha_1,\beta_1,\ldots,\alpha_g,\beta_g$ such that
$\langle \alpha_i,\alpha_j\rangle=\langle \beta_i,\beta_j\rangle=0$
and $\langle \alpha_i,\beta_j\rangle=\delta_{ij}$ for every $i,j \in
\{1,\ldots,g\}$, where $g=g(S)$ is the genus of $S$.

A submodule $A$ of $H_1 (S)$ is said to be \emph{Lagrangian} if the
intersection form of $S$ identically vanishes on $A \times A$.

%%%%%%%

\subsection{An obstruction to be weakly--Helmholtz} 
As usual, let $\Omega$ be a domain with smooth boundary and let
$S_0,\ldots,S_h$ be the connected components of $\partial
\Omega$. Since homology is additive with respect to the disjoint union
of topological spaces, we have a canonical isomorphism
$H_1(\partial\Omega) \cong \bigoplus_j H_1(S_j)$, which allows us to
define canonical projections $p_j\colon H_1 (\partial\Omega)\lra
H_1(S_j)$, $j=0,\ldots,h$. If $i_\ast\colon H_1 (\partial\Omega)\lra
H_1 (\overline\Omega)$ is the homomorphism induced by the inclusion,
we set
$$
P_j:=p_j ({\rm Ker} (i_\ast)) \subset  H_1(S_j),\qquad j=0,\ldots,h \, .
$$

\begin{lem}\label{obstruction} If $\Omega$ is weakly--Helmholtz, then 
$P_j$ is a Lagrangian submodule of $H_1(S_j)$ for every $j \in \{0,\ldots,h\}$.
\end{lem}
\Dim By Theorem~\ref{char1}, we can choose a basis of
$H_2(\overline\Omega,\partial\Omega)$ represented by a sy\-stem of
surfaces $\Ff=\{\Sigma_1,\ldots,\Sigma_r\}$. By Lemma~\ref{rA}, we
have that ${\rm Ker} (i_\ast)= {\rm Im}\, \partial$, where $\partial
\colon H_2(\overline\Omega,\partial\Omega)\lra H_1(\overline\Omega)$
is the usual ``boundary map'' of the sequence of the pair
$(\overline\Omega,\partial\Omega)$.  This readily implies that, for
every $j \in \{0,\ldots,h\}$, the module $P_j$ is generated by a set
of classes which are represented by \emph{pairwise disjoint}
$1$--cycles, whence the conclusion. \cvd

\begin{exa} \label{tubi}
{\rm As an application of the previous lemma, one can see that the
open tubular neighbourhood (homeomorphic to $S \times (0,1)$) of a
smooth surface $S$ of genus $g>0$ is {\it not} weakly--Helmholtz. In
fact, if $\gamma$ is any simple loop on $S\times\{1\}$, then the cycle
$(\gamma\times\{1\})\sqcup (- \gamma\times\{0\})$ bounds the annulus
$\gamma\times [0,1]$, so the class $[\gamma\times\{1\}]-
[\gamma\times\{0\}]$ lies in ${\rm Im}\,\partial={\rm Ker}
(i_\ast)$. After setting $S_i=S \times \{i\}$, $i=0,1$, we have then
$P_i=H_1(S_i)$, and $P_i$ is \emph{not} Lagrangian. In Figure
\ref{fig:storus}, it is drawn an open tubular neighbourhood of a
torus in $\R^3$ corresponding to the case $g=1$: such a domain is not
weakly--Helmholtz.}

\begin{figure}[htbp]
\begin{center}
 \includegraphics[height=5cm]{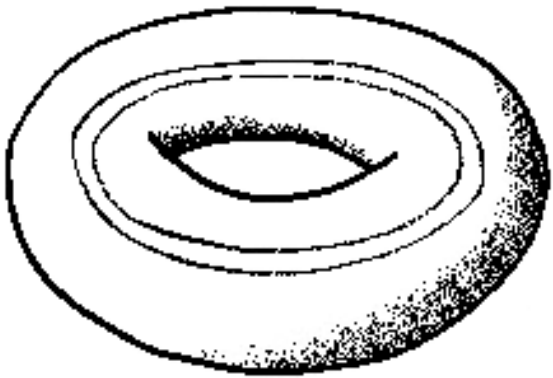}
\caption{\label{fig:storus} An open solid torus with a coaxial smaller
closed solid torus removed is not weakly--Helmholtz.}
\end{center}
\end{figure}
\end{exa}

The following lemma shows that, if $\partial\Omega$ is connected, then
Lemma~\ref{obstruction} does not provide any effective obstruction to
be weakly--Helmhlotz.

\begin{lem}\label{connected1} If the boundary $\partial \Omega = S_0$
is connected, then  ${\rm Ker}(i_*) \subset H_1(S_0)$ is a maximal Lagrangian
submodule of $H_1(S_0)$.
\end{lem} 
\Dim Lemma \ref{rA} implies that ${\rm Ker} (i_\ast)$ is a direct
summand of $H_1 (S_0)$ with ${\rm rank}\, {\rm
Ker}(i_*)=g(S_0)=\frac{{\rm rank}\, H_1(S_0)}{2}$, so it is enough to
show that ${\rm Ker} (i_\ast)$ is Lagrangian.

So let $\alpha$ be a $1$--cycle in ${\rm Ker} (i_\ast)$ represented by
a smooth loop $C_1\subset S_0$. If $[\beta]$ is any class in ${\rm
Ker} (i_\ast)={\rm Im}\, \partial$, then $[\beta]=\partial [\Sigma]$,
where $\Sigma$ is a properly embedded surface in
$(\overline\Omega,\partial\Omega)$. Since $S_0$ admits a collar in
$\overline\Omega$, we can push $\alpha$ a bit inside $\Omega$ and
obtain a $1$--cycle $\alpha'$ transverse to $\Sigma$. Since
$[\alpha']=i_\ast ([\alpha])=0$, the algebraic intersection between
$\alpha'$ and $\Sigma$ is null, and this easily implies in turn that
$\langle [\alpha],[\beta]\rangle =0$, whence the conclusion. \cvd

%%%%%%%
 
\smallskip

\subsection{Weakly--Helmholtz links}
We have seen in Corollary~\ref{knots} that all 
knots and all the box--domains of knots are weakly--Helmholtz. On the
other hand, if $L$ is the Hopf link (see Figure \ref{fig:hopf} below,
on the left), then $\compl(L)$ is diffeomorphic to an open tubular
neighbourhood of the standard torus in $\R^3$, so $\compl(L)$ is not
weakly--Helmholtz (see Example~\ref{tubi}). The same is true for
$\B(L)$ (see Lemma~\ref{equivalence}). Lemma~\ref{linked} below
generalizes this result to a large class of links. We say that two
components $K_1$ and $K_2$ of $L$ are {\it algebraically unlinked} if
$K_1$ is homologically trivial in $\compl(K_2)$. It turns out that
$[K_1]=0$ in $H_1 (\compl(K_2))$ if and only if $[K_2]=0$ in
$H_1(\compl(K_1))$, so the definition just given is indeed symmetric
in $K_1$ and $K_2$.  Equivalently, $K_1$ and $K_2$ are algebraically
unlinked if and only if their {\it linking number} vanishes; moreover
the linking number can be easily computed by using any planar link
diagram as half the sum of the signs at the crossing points betweem the
two components (for all this matter see e.g.~\cite{rolfsen} Section D
of Chapter 5).  Clearly, if two components of $L$ are geometrically
unlinked (see Subsection \ref{link}), a fortiori they are also
algebraically unlinked. The Whitehead link (see Figure~\ref{fig:white}
above on the left) is a celebrated example with two components that
are algebraically, but not geometrically, unlinked. The components
$K_1$ and $K_2$ are said to be {\it algebraically linked} it they are
not algebraically unlinked. Evidently, the Hopf link has algebraically
linked components.

\begin{lem}\label{linked}
If $L$ has algebraically linked components, then it is not
weakly--Helmholtz.
\end{lem}
\Dim Take two algebraically linked components $C_0$ and $C_1$ of $L$
and let $F_0$ be an oriented Seifert surface for $C_0$. As usual, we
can assume that $F_0$ is transverse to $C_1$ and to the corresponding
toric boundary component $S_1$ of $\partial \compl(L)=\partial U(L)$,
where $\compl(L)=S^3 \setminus U(L)$. Then the class $[\alpha]=p_1
(\partial [F_0 \setminus {\rm Int}(U(L))])\in p_1 ({\rm Ker}
(i_\ast))\subset H_1 (S_1)$ is represented by the oriented
intersection between $F_0$ and $S_1$, which is given by a finite
number of (possibly non--equioriented) copies of the meridian of
$S_1$. Since $C_0$ and $C_1$ are linked, the class $[\alpha]$ is not
null in $H_1 (S_1)$, and is therefore equal to a non--trivial multiple
of the class represented by the meridian of $S_1$. On the other hand,
also the class $[\beta]$ of the preferred longitude on $S_1$,
determined by any Seifert surface of $C_1$, belongs to $p_1 ({\rm Ker}
(i_\ast))$, and $\langle [\alpha],[\beta]\rangle \neq 0$, so Lemma
\ref{obstruction} implies that $L$ is not weakly--Helmholtz.  \cvd

\smallskip
%qui
\begin{remark}
{\rm Lemma~\ref{linked}
implies the Hopf link is not weakly--Helmholtz.
Thanks to
the Lemma~\ref{equivalence}, it follows that the box--domain of such a link,
drawn in the Figure \ref{fig:hopf} (on the right), is not
weakly--Helmholtz as well.}

%\vspace{-1em}

\begin{figure}[htbp]
\begin{center}
 \includegraphics[height=4.6cm]{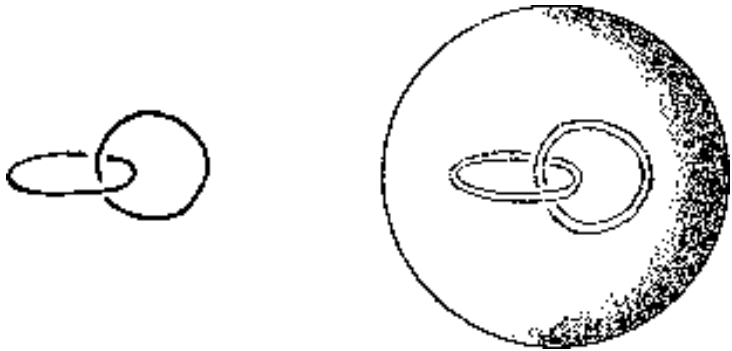}%{hopf-new.eps}
% \vspace{-1em}
\caption{\label{fig:hopf} {\rm A box--domain of a Hopf link is not
weakly--Helmholtz.}}
\end{center}
\end{figure}
\end{remark}

\smallskip

The following lemma considers the case of links with unlinked components.

\begin{lem}\label{unlinked} 
Suppose that the components $C_0,\ldots,C_k$ of a link $L$ are
algebraically unlinked with each other. Then there exists a family of
properly embedded surfaces $F_0,\ldots,F_k$ such that each $F_j$ is a
Seifert surface for $C_j$ and, if $i \neq j$, then $F_i$ and $F_j$
$($transversely$)$ intersect only in $\compl(L)$. Moreover, if
$i_j\colon S_j\lra \overline{\compl(L)}$ is the inclusion of the
boundary component corresponding to $C_j$ and $Q_j$ is the kernel of
$(i_{j})_\ast\colon H_1(S_j)\lra H_1(\overline{\compl(L)})$, then
$Q_j$ is generated by $($the class of$\, )$ the preferred longitude of
$C_j$, and ${\rm Ker} (i_*) = \bigoplus _j Q_j$.
\end{lem}
\Dim Fix $j\in \{0,\ldots,k\}$ and take an arbitrary Seifert surface
$F'_j$ of $C_j$ transverse to eve\-ry $C_h$, $h\neq j$. Up to
re--defining $\compl(L)$ as the complement in $S^3$ of smaller tubular
neighbourhoods of the $C_h$'s, we may also assume that, for each fixed
$h \neq j$, $F'_j$ intersects transversely each $S_h$ in a finite
number $m_1,\ldots,m_l$ of (possibly non--equioriented) copies of the
meridian of $S_h$, in such a way that each $m_i$ bounds a $2$--disk
$D_i$ in the interior of $F'_j$. Since the algebraic intersection of
$C_j$ and $C_h$ is null, we also have $[m_1]+\ldots+[m_l]=0$ in $H_1
(S_h)$, so the number of positively oriented meridians occurring in
the oriented intersection $F'_j\cap S_h$ equals the number of
negatively oriented meridians in the same intersection.

Let us now remove the $D_i$'s, $i=1,\ldots,l$, from the interior of
$F'_j$. In this way, we obtain a properly embedded surface with more
boundary components. We can now glue in pairs the added boundary
components by attaching $l/2$ disjoint annuli parallel to $S_h$ to
$l/2$ pairs of meridians in $F'_j\cap S_h$ having opposite
orientations. After applying the procedure just described to every
$h\neq j$, we obtain the desired Seifert surface $F_j$ that misses all
the $S_h$, $h\neq j$.

Now, if $[l_j]\in H_1 (S_j)$ is the class of the preferred longitude
of $C_j$, then $[l_j]=\partial [F_j]$, so $[l_j]$ lies in $Q_j$ and
hence ${\rm rank}\, \bigoplus_j Q_j =k+1={\rm rank}\, {\rm Ker}
(i_\ast)$. Now the conclusion follows from the fact that $\bigoplus_j
Q_j$ is a full submodule of $H_1 (\partial\Omega)$. \cvd

\smallskip

We may wonder if the Seifert surfaces of the previous lemma can be
chosen to be pairwise disjoint.  A classical definition is in order (see
\cite[p. 137]{rolfsen}).

\begin{defi}\label{b-link}
{\rm A link $L$ is a {\it boundary link} if it admits a system of {\it
disjoint} Seifert surfaces of its components.}
\end{defi}

\begin{figure}[htbp]
\begin{center}
 \includegraphics[height=4cm]{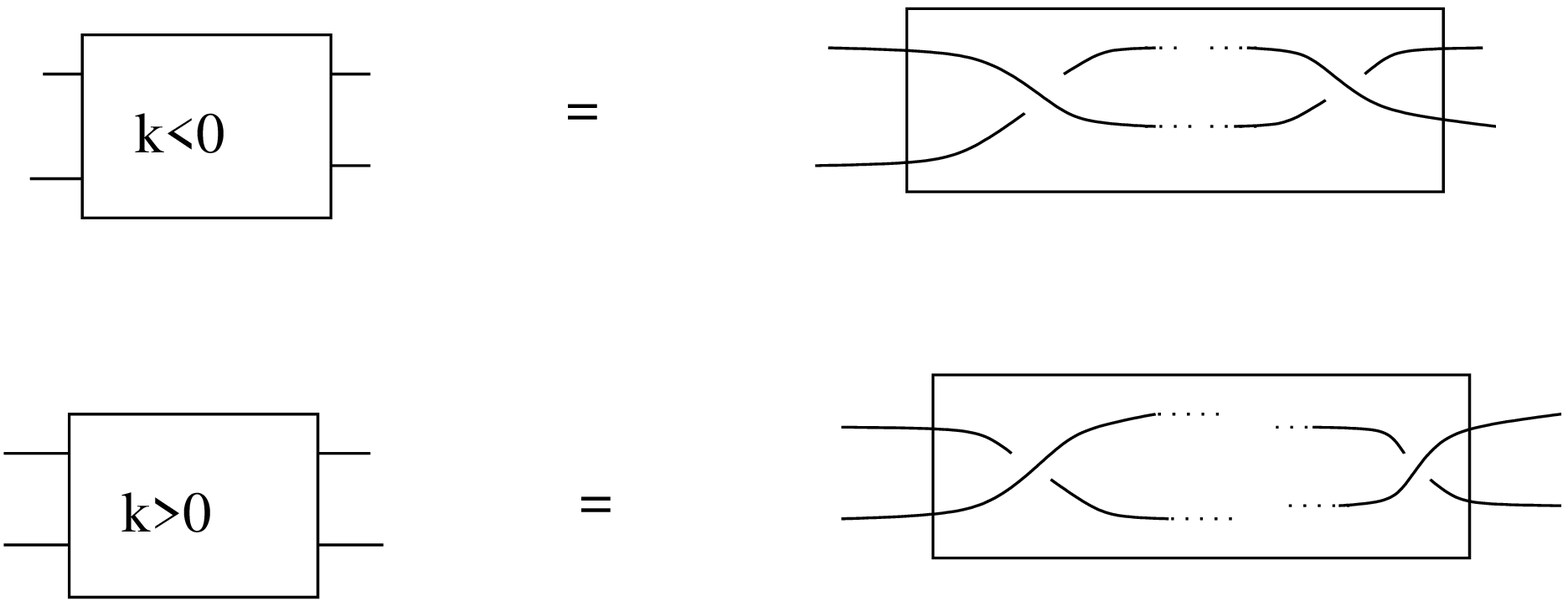}
\vspace{-1em}
\caption{\label{k-box} k--box.}
\end{center}
\end{figure}

\begin{figure}[htbp]
\begin{center}
 \includegraphics[height=6cm]{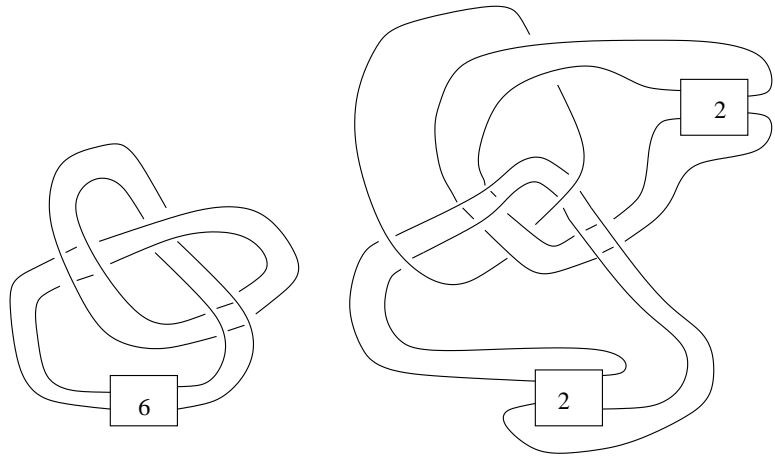}
\vspace{-1em}
\caption{\label{complicated_bl} Boundary links.}
\end{center}
\end{figure}

Of course, every knot is a boundary link. Every link $L$ with
geometrically unlinked components is a boundary link as well, as for every
component $C$, we can construct a Seifert surface contained in the
3--disk that separates $C$ from the other components
(see~\cite{rolfsen}).  However, there are boundary links that have 
geometrically linked components. For example every 2--components
links made by a non--trivial knot and its preferred longitude (recall
Example \ref{Seifert}) is a boundary link. On the left of Figure
\ref{complicated_bl}, we show the case of the trefoil knot, on the
right another more complicated 3--components boundary link (see
\cite{rolfsen} for other examples).  The meaning of the useful
square--boxes labelled by any integer $k$ is fixed in Figure \ref{k-box},
where it is understood that the box contains $|k|$ crossings.

On the other hand, the Whitehead link provides an example of a link with
algebraically unlinked components which is \emph{not} a boundary link (see
again~\cite[p. 137]{rolfsen}, and Example \ref{helm-nonhelm} below for an
even stronger result). So, in general, it is not possible to remove the
internal intersections of the Seifert surfaces provided by
Lemma~\ref{unlinked} by any local ``cut and paste'' procedure around the
intersection lines.

Let us now rephrase Theorem~\ref{char1} in the case of links.

\begin{cor}\label{char2-L}
A link $L$ with $r$ components is weakly--Helmholtz if and only if there
is a surjective homomorphism from $\pi_1(\Omega (L))$ to $\fr$.
\end{cor}
%\cvd

\smallskip

We recognize that the condition described in the last corollary is just
one current {\it definition} of {\it homology boundary
links}, so a link is weakly--Helmholtz if and only if it is a homology
boundary link. More precisely, putting together 
Corollary~\ref{char2-L} and Lemma~\ref{equivalence},
we obtain the following:

\begin{cor}\label{hb-L}
Given a link $L$ in $S^3$, the following assertions are equivalent:
\begin{itemize}
 \item[$(1)$] $L$ is weakly--Helmholtz.
 \item[$(2)$] $L$ is a homology boundary link.
 \item[$(3)$] $\B(L)$ is weakly--Helmholtz.
\end{itemize}
\end{cor}

Every classical boundary link is a homology boundary link. In fact, $L$ is
a boundary link if and only if there exists a surjective homomorphism
$\phi:\pi_1({\Omega (L)}) \lra \fr$, which furthermore (up to conjugacy)
sends the link meridians onto a set of generators of $\fr$. This
characterization of boundary links was originally given in~\cite{Smythe}
(see also \cite{Gutier}), where also  the relaxed definition of homology
boundary links was introduced.

Homology boundary links are an intriguing, very important class of links
widely studied in Knot Theory (the interested reader can find more
It is a nice occurence that our
discussion originated  from Helmholtz cuts, eventually leads to such a
distinguished class of links.

%%%

\subsection {On general weakly--Helmholtz domains}\label{general}
Getting an exhaustive description of weakly--Helmholtz domains,
similar to the characterization of Helmholtz ones
given in Theorem~\ref{helm-char}, looks somehow hopeless. This already
holds true for the special case of links. Note that concretely given a
link $L$ (for instance by means of a usual planar link diagram), with
algebraically unlinked components, it is in general a quite hard task
to decide whether or not it is homology boundary (for example, 
some non--trivial argument is needed even for showing that
the Whitehead link is not weakly--Helmholtz
-- see the examples below). The general case is even more complicated.
Up to ``Fox reimbedding'' (see Theorem~\ref{Fox}), 
it is not restrictive to deal with
domains $\Omega$ that are the complements of links of handlebodies
considered up to isotopy. As every handlebody is the regular
neighbourhood of a spine, which is a compact graph embedded in $S^3$
(i.e. a {\it spatial graph}), if $\Gamma$ is a {\it link of spines},
then we can naturally extend our previous notation and denote by $\compl
(\Gamma)$ the complement--domain of $\Gamma$. In the
case of a classical link $L$, i.e. in the case of a link of genus 1
handlebodies, we have in some sense a ``canonical'' spine
for $\compl (L)$: the link $L$ itself. This is no
longer true in the general case, in the sense that a link of
handlebodies, considered up to isotopy, can admit essentially
different links of spines. This represents a further complication in
the study of general weakly--Helmholtz domains.

\begin{figure}[htbp]
\begin{center}
 \includegraphics[height=1cm]{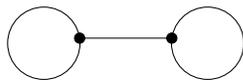}
\vspace{-1em}
\caption{\label{handcuff} Planar handcuff graph.}
\end{center}
\end{figure}

To illustrate the last claim, we will consider the simplest case of just
one genus 2 handlebody $\overline H$. Every such handlebody admits a spine
$\Gamma$, which is a spatial embedding of the so--called ``handcuff graph'' (a
planar realization of which is shown in Figure \ref{handcuff}).

If we remove from  $\Gamma$ the interior of the edge that connects the two
cycles (i.e. the ``isthmus'' of  $\Gamma$), then we get a classical link
$L_\Gamma$ with two components. Set $\Omega = \compl
(\Gamma)$ and  $ \Omega' =  \compl (L_{\Gamma})$. Clearly
$\overline{\Omega} \subset \overline{\Omega'}$, as the first is obtained by
removing a 1--handle from the second. 

The following proposition will allow to contruct many examples of
both non--homology boundary links with two algebraically unlinked
components, and  knotted genus 2 handlebodies having weakly--Helmholtz 
complementary domain.

\begin{prop}\label{link-link} 
With the notations just introduced, the following results hold:
\begin{enumerate}
 \item 
If  $L_\Gamma$  is a homology boundary
link, then $\Omega$ is weakly--Helmholtz.
\item 
Suppose that $\overline H$ is unknotted. Then $L_\Gamma$ is a homology boundary
link if and only if  $\Gamma$ is planar. In particular, 
if  $L_\Gamma$ is non--trivial, then it is
not a homology boundary link.
\end{enumerate}
\end{prop}

\Dim By a general position argument, it is easy to see that every loop
in $S^3\setminus L_\Gamma$ is homotopic to a loop that does not
intersect the isthmus of $\Gamma$. This implies that $i_*:
\pi_1(\overline{\Omega}) \to \pi_1(\overline{\Omega'})$ is
surjective. Then (1) follows immediately from Theorem \ref{char1} and
Corollary \ref{char2-L}.

Let us now suppose that $\overline H$ is unknotted. Then also
$\overline{\Omega}$ is an unknotted genus 2 handlebody, hence
$\pi_1(\overline{\Omega}) \cong \Z^{*2}$, and $\pi_1
(\overline{\Omega'})$ is isomorphic to a quotient of $\pi_1(\overline
\Omega )$. Therefore, if $L_\Gamma$ is homology boundary, then we have
a sequence of surjective homomorphisms
$$
\Z^{*2}\cong \pi_1(\overline{\Omega} )\to \pi_1 (\overline{\Omega'})\to
\Z^{*2}\ .
$$ But free groups are Hopfian (see~\cite{MKS}), which means that
every surjective homomorphism of $\Z^{*2}$ onto itself is in fact an
isomorphism, and this implies here that $\pi_1 (\Omega')$ is
isomorphic to $\Z^{*2}$. Under this hypothesis, a generalization to
links (see for instance Theorem 1.1 in \cite{hillman}) of
Papakyriakopoulos unknotting theorem for knots \cite{papa} (which is
based on his famous ``loop theorem'' -- see also \cite{rolfsen})
ensures that $L_\Gamma$ is trivial, and we can finally apply the
planarity results of \cite{scharlemann2} and conclude that $\Gamma$ is
planar.  \cvd
\smallskip

We stress that $\overline H$ may admit {\it infinitely many} handcuff
spines with pairwise {\it non--}isotopic associated links (see the
examples below). Hence, if $\Omega$ is not
weakly--Helmholtz, point (1) of the above proposition implies that
no such link is homology boundary.  However, checking whether this last
condition is satisfied seems to be very demanding.

\begin{exa}\label{helm-nonhelm}{ \rm
(1) In Figures \ref{handcuff1} and \ref{handcuff2}, we show some spatial
handcuff graphs $\Gamma$ that become planar via a finite sequence of spine
modifications that keep the handlebody  $\overline H$ fixed up to isotopy.
In Figure \ref{handcuff2}, it is understood that $h=(-1)^{k}2$.

\begin{figure}[htbp]
\begin{center}
 \includegraphics[height=5cm]{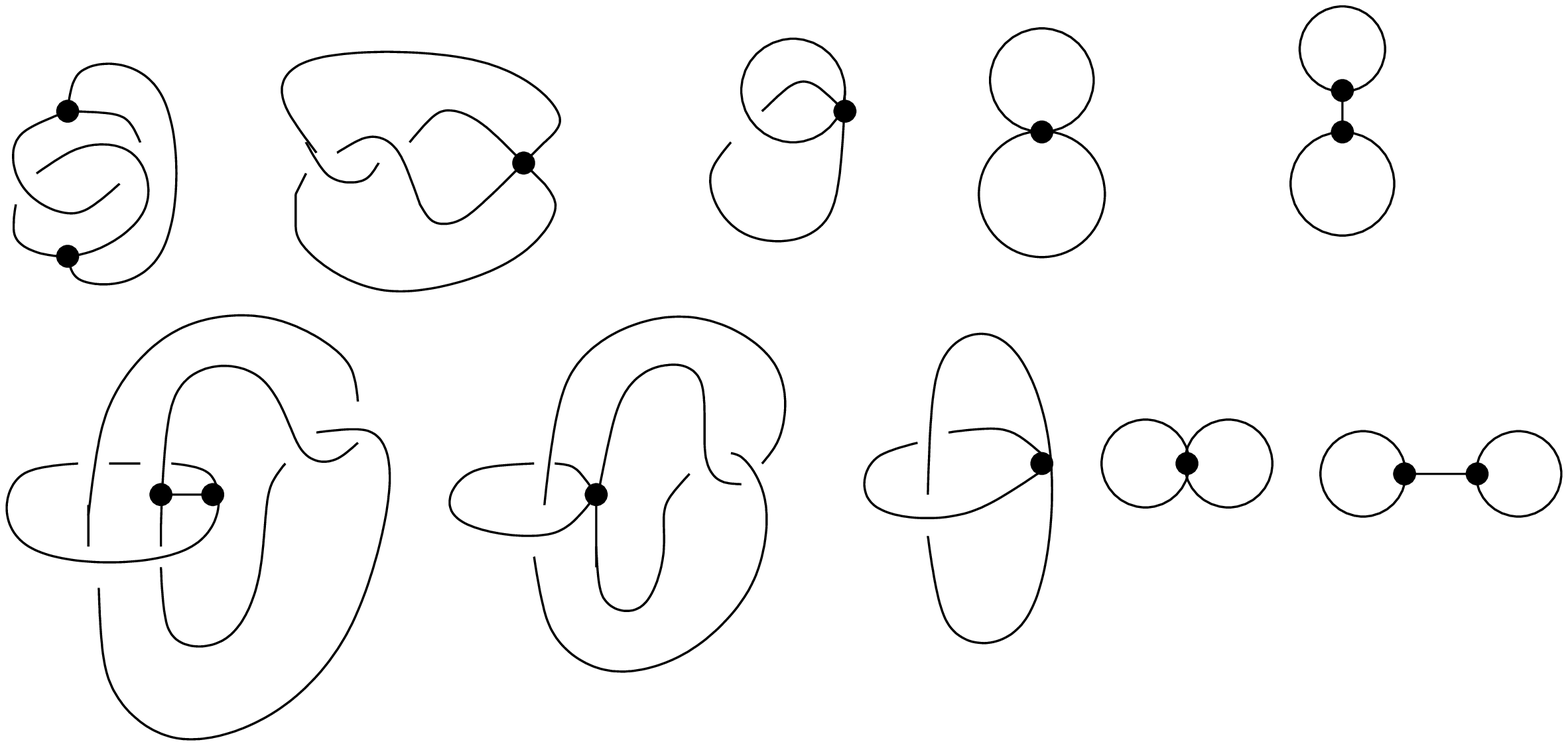}
\vspace{-1em}
\caption{\label{handcuff1} Unplanar vs planar handcuff spines.}
\end{center}
\end{figure}

\begin{figure}[htbp]
\begin{center}
 \includegraphics[height=5cm]{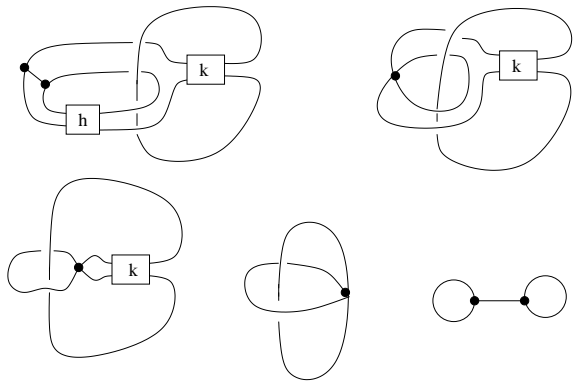}
\vspace{-1em}
\caption{\label{handcuff2} More unplanar vs planar handcuff spines.}
\end{center}
\end{figure}

The fact that the spines described here can be modified into planar
graphs shows that, in every case, $\overline H$ is unknotted, so, by point
(2) of Proposition~\ref{link-link}, we see that all the corresponding
non--trivial links $L_\Gamma$ are not homology boundary.  The first
example deals once again with the Hopf link by showing also the
somewhat non--intuitive phenomenon that being $L_\Gamma$ geometrically
linked does not prevent $\overline H$ to be unknotted. The second example
establishes that eventually the Whitehead link is not homology
boundary.  The examples described in Figure~\ref{handcuff2} provide an
infinite family of links (with the exceptions of $k=0, 1$ that produce
the trivial link) having algebraically unlinked components that are
not homology boundary. Note that every link in the family has one
unknotted component, while the other component is equal to the trefoil
knot when $k=-1$, the figure--eight knot when $k=-2$, etc. Note that when
$k=2$ we get 
Whitehead link again.
\smallskip

(2) If $L_\Gamma$ has geometrically unlinked components (i.e.~if it is
a split--link), then $\Omega = \compl (\Gamma)$ is weakly--Helmholtz
by point (1) of Proposition~\ref{link-link}. If we assume furthermore
that $L_\Gamma$ is non--trivial, then $\overline H$ is knotted by point
(2).  Remarkably, there exists also an example where $\overline H$ is
knotted whereas $L_\Gamma$ is trivial. In fact, it is proved in
\cite{ishii2} that the handlebody $\overline H$ determined by the spine
$\Gamma$ of Figure \ref{handcuff3} is knotted.

\begin{figure}[htbp]
\begin{center}
 \includegraphics[height=2cm]{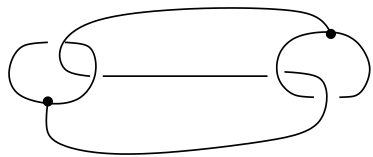}
\vspace{-1em}
\caption{\label{handcuff3} Knotted handlebody vs weakly--Helmholtz 
domain.}
\end{center}
\end{figure}

}
\end{exa}

Proposition~\ref{link-link} does not suggest how to construct examples  of
domains with connected boundary of genus 2 which are not weakly--Helmholtz. 
In fact we conclude our discussion
with the following open problem (as far as we know):

\begin{ques}{\rm Construct (if any) a 
knotted handlebody of genus 2 whose complement--domain
is not weakly--Helmholtz. Same question with arbitrary
genus. Due to Fox reimbedding Theorem, such handlebodies exist 
if and only if there exist domains with
{\it connected} locally flat boundary, which are not weakly--Helmholtz.}
\end{ques}

\subsection{Appendix}\label{appendix}
Without any pretension of being exhaustive, in this appendix we will
indicate to the interested reader some more advanced topics related to the
previous discussion.

Let $\Ll$ be a link of spines and suppose we are given a concrete presentation
of $\Ll$
(for instance by means of planar diagrams associated to generic planar
projections). Then it is rather easy to produce {\it finite presentations} of
the fundamental group of $S^3\setminus \Ll$ (such as the {\it Wirtinger
presentation} -- see \cite{rolfsen}). Fox's {\it free differential calculus}
\cite{Fox2} is a fundamental tool for the study of groups defined by
generators and relations. However, determining the corank starting from a
finite presentation of a group is in general a quite hard task. In
\cite{Stallings}, either this is done for certain presentations with
particular formal properties,  or one gives equivalent topological
3--dimensional reformulations, very close, in our framework, to the spirit
of Theorem \ref{char1}.

In the case of classical links, we may recur to certain, in principle
computable, increasingly discriminating sequence of invariants
(``obstructions'') whose vanishing is a necessary condition in order to be
homology boundary (at the initial step we have just the obstruction given
by the linking numbers of pairs of link components, discussed in
Lemma \ref{linked}). The original definition of such invariants is 
given in~\cite{MIL4}, so
that they are known as {\it Milnor's $\bar{\mu}$ invariants}. Let us
recall some of their formal features. For every integer $q>1$, for every
link $L$ with $N$ {\it ordered and oriented} components
$K_1,\dots , K_N$, for every $(l_1, \dots , l_p)\in \N^p$, with $1\leq
l_i \leq N$, $p<q$, it is defined an invariant of the form
$$ \bar{\mu}( l_1, \dots , l_p)(L)=[\mu( l_1, \dots , l_p)(L)]\in
\Z/\Delta( l_1, \dots , l_p)\Z$$
where:
\smallskip

- the integer $l_j$ is intended as a label of the component $K_{l_j}$
(note that any index $l_j$ can be repeated);
\smallskip

- the integer  $\mu( l_1, \dots , l_p)(L)$ is (not uniquely) obtained by
means of a determined procedure;

\smallskip

- the integer $\Delta( l_1, \dots , l_p)$ is defined inductively as the
g.c.d. of the numbers  $\mu( j_1, \dots , j_s)(L)$ where $s\geq 2$ and  $(
j_1, \dots , j_s)$ ranges over all cyclic permutations of proper
subsequences of   $( l_1, \dots , l_p)$.
\smallskip

- if $j_1\neq j_2$, the value
$\mu( j_1, j_2)(L)$ is the linking number
of the corresponding components.
\medskip

Strictly speaking, Milnor's invariants are isotopy invariants for ordered and oriented links. However, their vanishing does not depend on the chosen
 order or
orientation. The actual definition has a strong algebraic flavour, by
dealing with presentations of the fundamental group $G_1:=\pi_1(\compl
(L))$. Roughly speaking, Milnor's invariants detect
whether or not the (preferred)
longitudes of the link components can be expressed as longer and longer
commutators (i.e.~they detect how deep the longitudes live in
the {\it lower central series} of the link group, which is inductively defined
as follows: $G_1=\pi_1 (\compl (L))$, and 
$G_{n}=[G_{n-1},G_1]$ is the subgroup of $G$ generated by the set
$\{aba^{-1}b^{-1}\ ;\ a\in G_{n-1},\ b\in G_1\ \} $.
The invariants
relative to a given $q$ as above represent obstructions to the fact that the
longitudes belong to $G_q$).

In \cite{porter} and \cite{turaev}, it is established an equivalent
definition of Milnor's invariants in terms of the {\it Massey products} in
the systems $\{S^3\setminus K_{l_j}\}_{j=1}^p$. This approach provides an increasingly discriminating sequence of algebraic--topological obstructions defined by means of the cup product on singular 1--cochains with coefficients in
$\Z/\Delta( l_1, \dots , l_p)\Z$, and the coboundary operator.

In  \cite{Cochran0}, one can find a more geometric approach to these
invariants, based on the construction of so--called ``derived links''. This
method is particularly suited in order to deal with the ``first non--vanishing'' invariant (if any). In a sense it is a geometric realization of
the Massey products, working with relative 2--cycles rather than
1--cochains, and replacing the cup product with the transverse intersection
of such 2--cycles. The naive idea of a derived link is as follows. Consider
a link $L$ as in Lemma \ref{unlinked}, then we can construct a system of
Seifert surfaces intersecting transversely only in $\compl(L)$. We can manage in order that the intersection of each couple of surfaces is one connected knot in  $\compl (L)$. Each such knot splits in two parallel copies by slightly isotoping it out of both surfaces, by using the respective collars in the positive normal direction (accordingly to the orientations). By taking all knots obtained in this way, we get a derived link $L'$ of the given link $L$. We can define ``higher order'' invariants of $L$ by using the linking numbers of the pairs of components of $L'$. If all these linking numbers vanish, we iterate the procedure.

In \cite{Cochran0} and \cite{porter}, one finds some examples of
computations of non--trivial Milnor's invariants. In particular,
when $L$ is the Whitehead link, we see that $\bar{\mu}(1,1,2,2)(L)=1$, accordingly to the fact that $L$ is not homology boundary.

Milnor's invariants with pairwise distinct indices $l_j$ have a
particular meaning. In fact, they are invariant up to {\it link homotopy
equivalence} (also introduced in \cite{MIL4}). This means that one allows
homotopy with self--crossings of each link component, while crossings of
different components are not allowed. If a link $L$ is link homotopy
equivalent to a trivial link, then all such special Milnor's invariants
vanish. Note, for example, that the Whitehead link becomes trivial just by
performing one self--crossing at one component (see Figure~\ref{fig:white}). 
In \cite{fenn} or \cite{dimov}, it is proved that every boundary link is link homotopic to a trivial link.  The link--homotopy classification was given in \cite{MIL4} for 2-- and 3--components links, in \cite{levine} for 4--components ones; finally for all links in \cite{habbeger}. General Milnor's invariants are
invariants up to link {\it concordance} equivalence. Homology
boundary links have been widely studied in this framework 
\cite{Cochran0, Cochran, Cochran2, Cochran3, casson}.

The theory of links of spatial graphs as well as of links of handlebodies
is considerably less developed than the classical link theory. Links of spatial
graphs have been more intensively considered, by extending to them
different equivalence relations (``homotopy'', ``cobordism'', ``homology'',...)
\cite{fleming, taniyama, shinjo}. Particular efforts have been dedicated
to detect whether or not a link is planar (up to isotopy)
\cite{scharlemann1, Zhao}. A largely diffused approach consists in
associating to every link of spatial graphs some invariant families of
classical links \cite{kauffman, kauffman2, Ghuman}, in order to exploit
such a more developed theory.

The theory of links of handlebodies is even less developed. A natural
approach consists in considering links of spines, that is links of spatial
graphs up to isotopy coupled with suitable spine modifications that do not
alter the carried handlebodies \cite{scharlemann2, ishii1, ishii2}.

\bigskip

{\bf Note on the bibliography.} References \cite{ABDG} to \cite{temam}
form the ``Section A'' relative to Electromagnetism, Hydrodynamics and
Elasticity on domains in $\R^3$. References \cite{ALEX2} to
\cite{Zhao} form the ``Section B'' on (3--dimensional)
Differential/Algebraic/Geometric Topology.

\end{document}